\documentclass[letterpaper, paper,12pt]{AAS}	

\usepackage{bm}
\usepackage{amsmath}
\usepackage{subfigure}
\usepackage[colorlinks=true, pdfstartview=FitV, linkcolor=black, citecolor= black, urlcolor= black]{hyperref}
\usepackage{footnpag}	
\usepackage{placeins}
\usepackage{gensymb}

\usepackage{hyperref}
\hypersetup{
    colorlinks=true,
    linkcolor=blue,
    filecolor=magenta,      
    urlcolor=blue,
    pdftitle={Overleaf Example},
    pdfpagemode=FullScreen,
    }

\usepackage{algorithm}
\usepackage{algpseudocode}

\usepackage{graphicx}
\usepackage{multirow}

\PaperNumber{24-465}

\begin{document}


\title{Expanding the Class of Quadratic Control-Lyapunov Functions for Low-Thrust Trajectory Optimization}


\author{Nicholas P. Nurre\thanks{Graduate Student, Department of Aerospace Engineering, Auburn University. AAS Member.}, Saeid Tafazzol\thanks{Graduate Student, Department of Aerospace Engineering, Auburn University.}, Ehsan Taheri\thanks{Assistant Professor, Department of Aerospace Engineering, Auburn University. AAS Member.}}

\maketitle{}

\begin{abstract}
Control laws derived from Control-Lyapunov Functions (CLFs) offer an efficient way for generating near-optimal many-revolution low-thrust trajectories. A common approach to constructing CLFs is to consider the family of quadratic functions using a diagonal weighting matrix. In this paper, we explore the advantages of using a larger family of quadratic functions. More specifically, we consider positive-definite weighting matrices with non-zero off-diagonal elements (hereafter referred to as ``full'' matrices). We propose a novel eigendecomposition method for parameterizing \(N\)-dimensional weighting matrices that is easy to implement and guarantees positive-definiteness of the weighting matrices (\href{https://github.com/saeidtafazzol/positive_definite_parameterization}{companion code is provided}\footnote{\url{https://github.com/saeidtafazzol/positive_definite_parameterization}}). 
We use particle swarm optimization, which is a stochastic optimization algorithm, to optimize the parameters and generate near-optimal minimum-time low-thrust trajectories. Solutions obtained using a full positive-definite matrix are compared to the results from the (standard) diagonal weighting matrix for a number of benchmark problems. Results demonstrate that improvements in optimality are achieved, especially for maneuvers with large changes in orbital elements.
\end{abstract}

\section{Introduction}
Lyapunov control (LC) laws form an efficient class of methods for generating near-optimal (i.e., with respect to time and propellant consumption) many-revolution low-thrust trajectories. Some well-known control laws include Q-law \cite{petropoulos_low-thrust_2004} and the Chang-Chichka-Marsden (CCM) control law \cite{eui_chang_lyapunov-based_2002}. The key idea in developing LC laws is based on the Lyapunov stability theory in that, given a candidate Lyapunov function that is positive definite everywhere except at a prescribed ``goal'' (i.e., an equilibrium), the system will converge to the equilibrium if the Lyapunov function's time derivative is negative-definite. 

Typically, in low-thrust trajectory optimization applications, the candidate Lyapunov function, also known as the control-Lyapunov function (CLF), is defined as a positive-definite function of the errors between instantaneous and target orbital elements such that an orbit transfer maneuver is achieved when the CLF becomes zero. 

Let \(V\) define a candidate CLF, i.e., a function that is positive definite in the state space except at the system's equilibrium/origin (i.e., the desired target state). Quadratic CLFs are typically written as, 
\begin{equation} \label{eq:lyap}
    V = \frac{1}{2}\bm{w}^ \top\bm{K}\bm{w},
\end{equation}
where \(\bm{K}\) is a positive-definite parameter matrix and $\bm{w}$ denotes an error vector. Deriving a LC law from a candidate CLF essentially reduces to the task of finding the expression for the control vector that point-wise minimizes the time derivative of the CLF \cite{schaub_analytical_2018}. 

Ref.~\cite{atmaca_near-optimal_2024} provides an explanation for deriving the appropriate LC law and the stability of the resulting system in the context of low-thrust trajectory optimization problems. One of the main challenges in finding low-thrust transfer solutions with LC lies, mostly, in choosing a candidate CLF. A straightforward way to ensure \(V\) remains positive-definite is to consider a diagonal weighting matrix. This particular parameterization offers simplicity in that selecting positive elements along the diagonal entries of the weighting matrix automatically guarantees positive-definiteness of the weighting matrix. This is the structure that, for example, Q-law \cite{petropoulos_low-thrust_2004} and CCM control laws are based on. 

Q-law is the most popular low-thrust control law and many studies have focused on developing efficient methods to optimize the parameters of the Q-law to achieve either time- or fuel-optimal solutions. For example, Ref.~\cite{lee_design_2005} optimizes the parameters of the Q-law with two evolutionary/stochastic optimization algorithms (i.e., genetic algorithm and simulated annealing) and compares their performances. A mechanism for coasting is introduced in Ref.~\cite{petropoulos_low-thrust_2004} for Q-law, which allows for obtaining suboptimal minimum-fuel solutions. This coasting mechanism was leveraged in Ref.~\cite{lee_design_2005} to generate a Pareto front of time- and fuel-optimal solutions. Trajectories were obtained that were comparable to those obtained using other low-thrust trajectory optimization algorithms, such as NASA's Mystic. Ref.~\cite{holt_optimal_2021} uses reinforcement learning to determine non-constant state-dependent Q-law parameters that result in improved optimality as opposed to when constant parameters are used. Ref.~\cite{hecht_q-law_2024} extends Q-law to account for Sun-angle constraints in the trajectory design of solar electric propulsion, demonstrating the versatility of LC laws to handle practical operational constraints. Some studies, e.g., Refs.~\cite{alizadeh_static_2011} and \cite{atallah_inverse-optimal_2024}, have explored the feasibility of using solutions obtained through LC to initialize the costates in indirect methods by leveraging the connection between Lyapunov functions and the Hamilton-Jacobi-Bellman equation -- an idea formalized under the so-called inverse optimal control theory \cite{freeman_inverse_1996}.

LC remains a state-of-the-art method for rapid low-thrust trajectory design and optimization. However, to the best knowledge of the authors, no work explores using a full weighting matrix (i.e., one that considers cross-product terms in addition to perfect square terms in \(V\)) in the CLF. Considering the cross-product terms can enlarge the possible set of solutions and possibly result in improved solutions in terms of optimality. Such a method that allows for positive-definiteness on \(\bm{K}\) with simple bound constraints on the decision variables would allow for not only constant weighting matrices to be efficiently optimized over, but also allow for non-constant weighting matrices (like in Ref.~\cite{holt_optimal_2021}) to be efficiently found.

The main contributions of this paper are as follows: we explore enlarging the family of quadratic CLFs for low-thrust trajectory optimization by considering a full symmetric weighting matrix \(\bm{K}\) (see Eq.~\eqref{eq:lyap}). To efficiently search for symmetric positive-definite \(\bm{K}\) matrices, eigendecomposition \cite{strang_introduction_2023} is used to parameterize the weighting matrix \(\bm{K}\) through its eigenvalues and eigenvectors. Positive definiteness can be easily guaranteed through eigendecomposition by bounding the eigenvalues of the \(\bm{K}\) matrix to be positive. Because eigenvectors corresponding to real, unique eigenvalues are orthogonal, the matrix of eigenvectors characterizes a rotation matrix. In 3 dimensions, such a rotation matrix can be parameterized and constructed with 3 Euler angles. However, this must be generalized to higher dimensions to be applicable to a larger class of problems such as orbit transfer maneuvers in which up to five orbital elements are targeted. In this work, we propose a method for constructing the rotation matrices using a minimal-parameter representation. The proposed minimal-parameter representation allows for full positive-definite \(\bm{K}\) matrices of arbitrary size to be constructed and the entire parameter space to be searched over, enlarging the family of trajectories achievable with the derived control law. To our best knowledge, no existing work in low-thrust trajectory optimization considers this larger class of quadratic CLFs.

The remainder of the paper is organized as follows. First, the procedure for parameterization and construction of the symmetric positive-definite matrices are given. Then, this parametrization is applied to a set of benchmark orbit transfer problems \cite{petropoulos_low-thrust_2004}. Next, the results of these orbit transfers are presented when the two parameterizations of the $\bm{K}$ matrix, i.e., the standard diagonal and the full symmetric matrices denoted as, \(\bm{K}_1\) and \(\bm{K}_2\), respectively, are considered. Finally, a conclusion is given and some future works are outlined.

\section{Symmetric Positive-Definite Matrix Parameterization}

An effective parameterization that captures the entire space of positive-definite matrices is beneficial for LC. It is essential that such a parameterization does not introduce additional degrees of freedom than necessary. To achieve such a parameterization, we leverage the spectral theorem \cite{strang_introduction_2023} for parameterizing symmetric matrices, which can be written as,
\begin{equation} \label{eq:spectral}
    \bm{K} = \bm{Q} \bm{\Lambda} \bm{Q}^{\top}.
\end{equation}

In Eq.~\eqref{eq:spectral}, $\bm{K}$ is an $N \times N$ symmetric matrix. The matrix $\bm{Q}$ is orthonormal (also referred to as a rotation matrix) and the matrix $\bm{\Lambda}$ is a diagonal one. Such a representation of symmetric matrices is also referred to as eigendecomposition, where $\bm{\Lambda}$ contains the eigenvalues of $\bm{K}$ on its main diagonal. We can parameterize $\bm{\Lambda}$ by its diagonal entries, ensuring all eigenvalues are positive to maintain the positive definiteness of $\bm{K}$.

The matrix $\bm{Q}$, being a rotation matrix, can be parameterized in terms of angles. For instance, in a 2-dimensional case, $\bm{Q}$ can be parameterized by an angle $\theta$ as,
\begin{equation}
    \bm{Q} = \begin{bmatrix} \cos(\theta) & -\sin(\theta) \\
                        \sin(\theta) & \cos(\theta) \end{bmatrix}.
\end{equation}

In the 3-dimensional case, three Euler angles can be used to construct $\bm{Q}$. Inspired by these examples, a general $N$-dimensional method \cite{hoffman_generalization_1972} is used to construct the rotation matrix. To do so, we utilize the hyperspherical geometry to construct unit vectors using angles as,
\begin{equation}    
    \bm{v}_i^p = \begin{bmatrix}
        \cos(\theta_{i1}) \\ \sin(\theta_{i1}) \cos(\theta_{i2}) 
        \\ \cdots \\ \prod_{j=1}^{N-i-1} \sin(\theta_{ij}) \cos{(\theta_{i(N-i)})}\\
        \prod_{j=1}^{N-i} \sin(\theta_{ij})
    \end{bmatrix},
\label{eq:unit_vec}
\end{equation}
where $\bm{v}_i^p$ is a unit vector in an $N-i + 1$ dimensional space. The superscript $p$ is used to denote that  $\bm{v}_i$ is not expressed in the original $N$-dimensional space and it should be projected back to the original dimension to construct the desired rotation matrix. Since by definition $\bm{Q}$ is an orthonormal matrix, all of its constituent vectors should be mutually orthogonal, which equivalently means that the generated vectors should be in the null space of all previous vectors. 

Based on the mentioned fundamentals, Algorithm \ref{alg:rotation} is proposed to parameterize a rotation matrix $\bm{Q}$ using angles $\theta_{ij}$ for $i = 1 ,\ldots , N, \; j = 1, \ldots, N-i$ that will be used to construct all the symmetric positive-definite matrices in this paper. Minimal working MATLAB and Python code is provided \href{https://github.com/saeidtafazzol/positive_definite_parameterization}{here}\footnote{\url{https://github.com/saeidtafazzol/positive_definite_parameterization}}.

\begin{algorithm}
\caption{N-dimensional Rotation Matrix 
Parameterization}
\label{alg:rotation}
\begin{algorithmic}[1]
\Statex \textbf{Input}: Angles $\theta_{ij}$ for $i = 1 ,\ldots , N, \; j = 1, \ldots, N-i$ 
\Statex \textbf{Output:} an Orthonormal Matrix $\bm{Q}$
\Statex \textbf{Initialize:} $\bm{Q} \leftarrow \{\}$ \Comment{The set of orthonormal vectors}

\For{$i = 1$ \textbf{to} $N$}
    \State $S = \text{null-space}(\bm{Q}^{\top})$
    \State $B = \text{Basis} (S)$ \Comment{Gram–Schmidt process}

    \State Construct $\bm{v}_i^p$ using $\theta_{i \, (1, \ldots, N-i)}$ and Eq.~\eqref{eq:unit_vec}
    \State $\bm{v}_i = B \bm{v}_i^p$ \Comment{Project back to $N$-dimensional space}
    \State $\bm{Q} \leftarrow [\bm{Q},  \bm{v}_i]$ \Comment{Append $\bm{v}_i$ as a new column to $\bm{Q}$}
    
\EndFor

\end{algorithmic}
\end{algorithm}

The outlined procedure requires \(\frac{N(N-1)}{2}\) angles and \(N\) eigenvalues, resulting in a total number of \(\frac{N(N+1)}{2}\) parameters. This count matches the total number of parameters needed to construct a symmetric matrix, confirming that the procedure does not introduce any redundant degrees of freedom. Additionally, the only constraints that need to be enforced while optimizing the parameters are simple bound (box) constraints on the parameters.

\section{Derivation of the control law and definition of error terms}
To characterize the benefits of using a full weighting matrix in the CLF for generating low-thrust trajectories, a variety of low-thrust transfer problems are solved with two different control laws: one with a diagonal weighting matrix, \(\bm{K}_1\), and one with a full (symmetric) weighting matrix, \(\bm{K}_2\), which is generated using Algorithm \ref{alg:rotation}. It is assumed that the spacecraft is thrusting for the entirety of the transfer. The state of the spacecraft, \(\bm{x}\), consists of Cartesian position \(\bm{r}\) and velocity \(\bm{v}\)  vectors, defined in an inertial frame centered at the central body defined by the unit vectors \(\left\{\hat{\bm{x}},\,\hat{\bm{y}},\,\hat{\bm{z}}\right\}\). All acceleration vectors are expressed with respect to the inertial frame. An additional state, \(m\), is used to track the spacecraft's mass. The spacecraft is subject to the equations of motion defined as
\begin{align} \label{eq: eom}
    \bm{x} & = \begin{bmatrix} \bm{r} \\ \bm{v} \\ m \end{bmatrix}, & \dot{\bm{x}} & = \begin{bmatrix} \bm{v} \\ -\frac{\mu}{r^3} \bm{r} + \frac{T_\text{max}}{m}\hat{\bm{\alpha}} \\ -\frac{T_\text{max}}{I_\text{sp}g_0} \end{bmatrix},
\end{align}
where \(\mu\) is the gravitational parameter of the central body, \(r=\|\bm{r}\|\) is the spacecraft's distance from the central body, \(T_\text{max}\) is the spacecraft's maximum thrust, \(I_\text{sp}\) is the spacecraft's specific impulse, \(g_0\) is the acceleration of gravity on Earth at sea level defined as 9.80665 m/s\(^2\), and \(\hat{\bm{\alpha}}\) is the thrust unit steering vector assumed to freely orient in space but subject to the unit vector constraint, \(\hat{\bm{\alpha}}^\top\hat{\bm{\alpha}}=1\). 

In the considered trajectory design problems, the goal is to transfer the spacecraft starting from a fully defined state at time \(t_0\) to a state along an osculating orbit over the time horizon \(t\in[t_0,\infty)\) in minimum time. The boundary conditions and parameters for five transfer scenarios that we have solved are summarized in Table \ref{tab: orbit transfers}. These are a set of transfer cases originally defined in Ref.~\cite{petropoulos_low-thrust_2004} and have been solved as benchmark problems in a number of other studies, including Refs.~\cite{lee_design_2005} and \cite{varga_many-revolution_2016}. Many works consider eclipse-induced coasting and/or perturbations to the dynamics such as accelerations due to Earth's \(J_2\), but the present work does not for the sake of a simple comparison of the two control laws. However, our future work will include higher-fidelity models consisting of eclipses and perturbations. 

The error vector, \(\bm{w}\) (see Eq.~\eqref{eq:lyap}), for a candidate CLF should be a function of states that becomes $\bm{0}$ when the final boundary conditions in Table \ref{tab: orbit transfers} are satisfied. For Case E, the error vector is defined in terms of the specific angular momentum vector, \(\bm{h}\), and the eccentricity vector, \(\bm{e}\) (also known as the Laplace vector under a different scaling). Definitions of each are as follows \cite{vallado_fundamentals_2022}:
\begin{align} 
    \bm{h} & = \bm{r}\times\bm{v}, & \bm{e} & = \frac{\left(v^2-\frac{\mu}{r}\right)\bm{r} - \left(\bm{r}^\top\bm{v}\right)\bm{v}}{\mu}.
\end{align}

The error vector for Case E, \(\bm{w}_\text{E}\), is defined as, 
\begin{equation}
    \bm{w}_\text{E} = \begin{bmatrix} \bm{h} - \bm{h}_\text{T} \\ \bm{e} - \bm{e}_\text{T} \end{bmatrix},
\end{equation}
where the subscript `$\text{T}$' corresponds to the target orbit values. Note that the considered choice of the error vector variables is entirely arbitrary. In fact, a vector of only 5 dimensions is needed since only 5 orbital elements are being targeted at the final boundary condition. Other error vector choices may offer improved results, but the emphasis of this work is on the comparison of the two different weighting matrices, \(\bm{K}_1\) and \(\bm{K}_2\). 

For Cases A through D, the boundary conditions are defined in terms of a combination of a set of orbital elements consisting of the specific angular momentum magnitude \(h\), eccentricity \(e\), inclination \(i\), and right ascension of the ascending node \(\Omega\). Definitions of the orbital elements are as follows \cite{vallado_fundamentals_2022}:
\begin{align} 
    h & = \left\|\bm{h}\right\|, & e & = \left\|\bm{e}\right\|, & i & = \cos^{-1}{\left(\frac{h_{\hat{\bm{z}}}}{h}\right)}, & \Omega & = \cos^{-1}{\left(\frac{\hat{\bm{x}}\cdot\bm{n}}{\|\hat{\bm{x}}\|\|\bm{n}\|}\right)},
\end{align}
where \(\bm{n}=\hat{\bm{z}}\times\bm{h} = [n_x,n_y,n_z]^{\top}\) is the line of nodes vector and a quadrant check is performed on \(\Omega\) such that if \(n_{y}<0\) then \(\Omega = 360 \degree - \Omega\). 

The error vectors for Cases A through D are defined as follows:
\begin{align}
    \bm{w}_\text{A} & = \begin{bmatrix} h - h_\text{T} \\ e - e_\text{T} \end{bmatrix}, & \bm{w}_\text{B} & = \begin{bmatrix} h - h_\text{T} \\ e - e_\text{T} \\ i - i_\text{T} \end{bmatrix}, & \bm{w}_\text{C} & = \begin{bmatrix} h - h_\text{T} \\ e - e_\text{T}\end{bmatrix}, & \bm{w}_\text{D} & = \begin{bmatrix} h - h_\text{T} \\ e - e_\text{T} \\ i - i_\text{T} \\ \Omega - \Omega_\text{T} \end{bmatrix}.
\end{align}

Note again that these error vector choices are arbitrary. For instance, instead of \(\Omega - \Omega_\text{T}\), Ref.~\cite{petropoulos_low-thrust_2004} suggests using \(\cos^{-1}{\left(\cos{\left(\Omega - \Omega_\text{T}\right)}\right)}\) to provide the angular measure the ``short way round'' the circle the angle is being measured on.

\begin{table}[]
\centering
\caption{Orbit Transfers}
\label{tab: orbit transfers}
\resizebox{\textwidth}{!}{%
\begin{tabular}{|c||l|r|r|r|r|r|r|r|r|r|r|r|}
\hline
\multirow{3}{*}{Case} &
  \multicolumn{1}{c|}{\multirow{3}{*}{Orbit}} &
  \multirow{3}{*}{\begin{tabular}[c]{@{}r@{}}\(a\)\\ {[}kg{]}\end{tabular}} &
  \multirow{3}{*}{\begin{tabular}[c]{@{}r@{}}\(e\)\\ {[}-{]}\end{tabular}} &
  \multirow{3}{*}{\begin{tabular}[c]{@{}r@{}}\(i\)\\ {[}deg{]}\end{tabular}} &
  \multirow{3}{*}{\begin{tabular}[c]{@{}r@{}}\(\omega\)\\ {[}deg{]}\end{tabular}} &
  \multirow{3}{*}{\begin{tabular}[c]{@{}r@{}}\(\Omega\)\\ {[}deg{]}\end{tabular}} &
  \multirow{3}{*}{\begin{tabular}[c]{@{}r@{}}Thrust\\ {[}N{]}\end{tabular}} &
  \multirow{3}{*}{\begin{tabular}[c]{@{}r@{}}Initial\\ Mass\\ {[}kg{]}\end{tabular}} &
  \multirow{3}{*}{\begin{tabular}[c]{@{}r@{}}Specific\\ Impulse\\ {[}s{]}\end{tabular}} &
  \multirow{3}{*}{\begin{tabular}[c]{@{}r@{}}Central\\ Body\end{tabular}} &
  \multirow{3}{*}{\begin{tabular}[c]{@{}r@{}}PSO\\ Swarm\\ Size\end{tabular}} &
  \multirow{3}{*}{\begin{tabular}[c]{@{}r@{}}PSO\\ Max.\\ Iterations\end{tabular}} \\
 &
  \multicolumn{1}{c|}{} &
   &
   &
   &
   &
   &
   &
   &
   &
   &
   &
   \\
 &
  \multicolumn{1}{c|}{} &
   &
   &
   &
   &
   &
   &
   &
   &
   &
   &
   \\ \hline \hline
\multirow{2}{*}{A} &
  Initial &
  7000 &
  0.01 &
  0.05 &
  0 &
  0 &
  \multirow{2}{*}{1} &
  \multirow{2}{*}{300} &
  \multirow{2}{*}{3100} &
  \multirow{2}{*}{Earth} &
  \multirow{2}{*}{50} &
  \multirow{2}{*}{50} \\ \cline{2-7}
 &
  Target &
  42000 &
  0.01 &
  free &
  free &
  free &
   &
   &
   &
   &
   &
   \\ \hline
\multirow{2}{*}{B} &
  Initial &
  24505.9 &
  0.725 &
  7.05 &
  0 &
  0 &
  \multirow{2}{*}{0.350} &
  \multirow{2}{*}{2000} &
  \multirow{2}{*}{2000} &
  \multirow{2}{*}{Earth} &
  \multirow{2}{*}{50} &
  \multirow{2}{*}{50} \\ \cline{2-7}
 &
  Target &
  42165 &
  0.001 &
  0.05 &
  free &
  free &
   &
   &
   &
   &
   &
   \\ \hline
\multirow{2}{*}{C} &
  Initial &
  9222.7 &
  0.2 &
  0.573 &
  0 &
  0 &
  \multirow{2}{*}{9.3} &
  \multirow{2}{*}{300} &
  \multirow{2}{*}{3100} &
  \multirow{2}{*}{Earth} &
  \multirow{2}{*}{50} &
  \multirow{2}{*}{50} \\ \cline{2-7}
 &
  Target &
  30000 &
  0.7 &
  free &
  free &
  free &
   &
   &
   &
   &
   &
   \\ \hline
\multirow{2}{*}{D} &
  Initial &
  944.64 &
  0.015 &
  90.06 &
  156.90 &
  -24.60 &
  \multirow{2}{*}{0.045} &
  \multirow{2}{*}{950} &
  \multirow{2}{*}{3045} &
  \multirow{2}{*}{Vesta} &
  \multirow{2}{*}{50} &
  \multirow{2}{*}{50} \\ \cline{2-7}
 &
  Target &
  401.72 &
  0.012 &
  90.01 &
  free &
  -40.73 &
   &
   &
   &
   &
   &
   \\ \hline
\multirow{2}{*}{E} &
  Initial &
  24505.9 &
  0.725 &
  0.06 &
  0 &
  0 &
  \multirow{2}{*}{2} &
  \multirow{2}{*}{2000} &
  \multirow{2}{*}{2000} &
  \multirow{2}{*}{Earth} &
  \multirow{2}{*}{100} &
  \multirow{2}{*}{300} \\ \cline{2-7}
 &
  Target &
  26500 &
  0.7 &
  116 &
  270 &
  180 &
   &
   &
   &
   &
   &
   \\ \hline
\end{tabular}%
}
\end{table}

The CLF for each of the transfer problems is defined as
\begin{equation} \label{eq: CLF}
    V_{j,i} = \frac{1}{2}\bm{w}_j^\top\bm{K}_i\bm{w}_j,
\end{equation}
for \(i=\{1,2\}\) and \(j=\{\text{A},\text{B},\text{C},\text{D},\text{E}\}\). The control law is derived to make the total time derivative of Eq.~\eqref{eq: CLF} (note that indices \(j\) and \(i\) are not written for brevity),
\begin{equation}
    \dot{V} = \frac{\partial V}{\partial \bm{r}}\frac{\partial \bm{r}}{\partial t} + \frac{\partial V}{\partial \bm{v}}\frac{\partial \bm{v}}{\partial t} + \frac{\partial V}{\partial m}\frac{\partial m}{\partial t} + \frac{\partial V}{\partial t}=\frac{\partial V}{\partial \bm{r}}\bm{v}+\frac{\partial V}{\partial \bm{v}}\left(-\frac{\mu}{r^3} \bm{r} + \frac{T_\text{max}}{m}\hat{\bm{\alpha}}\right),
\end{equation}
as negative as possible (pointwise), leading to the control law, given in Eq.~\eqref{eq: control law}, as,
\begin{equation} \label{eq: control law}
    \hat{\bm{\alpha}}^*=-\left(\frac{\partial V}{\partial \bm{v}}\right)^\top/\left\|\frac{\partial V}{\partial \bm{v}}\right\|.
\end{equation}

The expression for thrust steering unit vector, given in Eq.~\eqref{eq: control law}, is derived using the automatic differentiation capabilities of CasADi \cite{andersson_casadi_2019}.


\section{Results}
Earth's gravitational parameter is 398600.49 km\(^3\)/s\(^2\) and Vesta's gravitational parameter is 17.8 km\(^3\)/s\(^2\). A canonical scaling is performed on all Cartesian states, such that \(\mu\) is 1. For the Earth-centered problems the distance unit (DU) is 6378.1366 km and the time unit (TU) is 806.8110 seconds and for the Vesta-centered problem, DU is 289 km and TU is 1164.4927 seconds.

The diagonal weighting matrix, \(\bm{K}_1\), is trivially constructed and the full weighting matrix \(\bm{K}_2\) is constructed following the steps outlined in Algorithm \ref{alg:rotation}. Particle swarm optimization (PSO) is used to optimize the parameters for each of the control laws. Because PSO is a stochastic optimization algorithm, it is ran 5 separate times for each control law. The swarm sizes and the number of iterations are summarized in Table \ref{tab: orbit transfers}. The eigenvalue parameters are arbitrarily bounded between 0 and 100. 
The simulations are run in MATLAB and the variable-step variable-order explicit differential equation solver \verb|ode113| is used to integrate Eq.~\eqref{eq: eom} with the control law in Eq.~\eqref{eq: control law} substituted in for \(\hat{\bm{\alpha}}\).

Because LC laws will not make the system converge in finite time, a tolerance \(\epsilon=1.0\times10^{-4}\) is defined such that all elements of \(\bm{w}\) being less than \(\epsilon\) corresponds to a converged solution. Note that this value is chosen arbitrarily so that an accurate orbit insertion is achieved, but other values can be chosen that provide a better balance between computational efficiency and orbit insertion accuracy. The equations of motion are integrated for an arbitrarily large time horizon, and the event detection capability of \verb|ode113| is used to determine when the insertion criterion (i.e., $\| \bm{w}\|_{\infty} \leq \epsilon$) is satisfied.

\begin{table}[htbp!]
\centering
\caption{Comparison of the results using the two parameter matrix formulations. Bold faced values denote the best solution for that particular transfer case and parameter matrix.}
\label{tab: results}
\resizebox{\textwidth}{!}{%
\begin{tabular}{|c||c|rrrrrr|}
\hline
\multirow{2}{*}{Case} &
  \multirow{2}{*}{\begin{tabular}[c]{@{}c@{}}Parameter\\ Matrix\end{tabular}} &
  \multicolumn{6}{c|}{Time of Flight {[}days{]}} \\ \cline{3-8} 
 &
   &
  \multicolumn{1}{c|}{Run 1} &
  \multicolumn{1}{c|}{Run 2} &
  \multicolumn{1}{c|}{Run 3} &
  \multicolumn{1}{c|}{Run 4} &
  \multicolumn{1}{c|}{Run 5} &
  \multicolumn{1}{c|}{Average} \\ \hline \hline
\multirow{2}{*}{A} &
  \(\bm{K}_1\) &
  \multicolumn{1}{r|}{14.5705} &
  \multicolumn{1}{r|}{14.5702} &
  \multicolumn{1}{r|}{\textbf{14.5700}} &
  \multicolumn{1}{r|}{14.5703} &
  \multicolumn{1}{r|}{14.5701} &
  14.5702 \\ \cline{2-8} 
 &
  \(\bm{K}_2\) &
  \multicolumn{1}{r|}{14.4751} &
  \multicolumn{1}{r|}{14.4748} &
  \multicolumn{1}{r|}{14.4749} &
  \multicolumn{1}{r|}{14.4754} &
  \multicolumn{1}{r|}{\textbf{14.4748}} &
  14.4750 \\ \hline
\multirow{2}{*}{B} &
  \(\bm{K}_1\) &
  \multicolumn{1}{r|}{\textbf{142.2285}} &
  \multicolumn{1}{r|}{142.2286} &
  \multicolumn{1}{r|}{142.2286} &
  \multicolumn{1}{r|}{142.2286} &
  \multicolumn{1}{r|}{142.2286} &
  142.2286 \\ \cline{2-8} 
 &
  \(\bm{K}_2\) &
  \multicolumn{1}{r|}{139.0297} &
  \multicolumn{1}{r|}{139.0674} &
  \multicolumn{1}{r|}{\textbf{139.0203}} &
  \multicolumn{1}{r|}{139.0739} &
  \multicolumn{1}{r|}{139.0235} &
  139.0430 \\ \hline
\multirow{2}{*}{C} &
  \(\bm{K}_1\) &
  \multicolumn{1}{r|}{1.5102} &
  \multicolumn{1}{r|}{1.5102} &
  \multicolumn{1}{r|}{1.5102} &
  \multicolumn{1}{r|}{1.5102} &
  \multicolumn{1}{r|}{\textbf{1.5102}} &
  1.5102 \\ \cline{2-8} 
 &
  \(\bm{K}_2\) &
  \multicolumn{1}{r|}{1.4918} &
  \multicolumn{1}{r|}{1.4918} &
  \multicolumn{1}{r|}{1.4918} &
  \multicolumn{1}{r|}{1.4918} &
  \multicolumn{1}{r|}{\textbf{1.4918}} &
  1.4918 \\ \hline
\multirow{2}{*}{D} &
  \(\bm{K}_1\) &
  \multicolumn{1}{r|}{24.9903} &
  \multicolumn{1}{r|}{24.9903} &
  \multicolumn{1}{r|}{24.9903} &
  \multicolumn{1}{r|}{24.9903} &
  \multicolumn{1}{r|}{\textbf{24.9903}} &
  24.9903 \\ \cline{2-8} 
 &
  \(\bm{K}_2\) &
  \multicolumn{1}{r|}{\textbf{24.6992}} &
  \multicolumn{1}{r|}{24.7059} &
  \multicolumn{1}{r|}{24.9456} &
  \multicolumn{1}{r|}{24.9456} &
  \multicolumn{1}{r|}{24.8868} &
  24.8366 \\ \hline
\multirow{2}{*}{E} &
  \(\bm{K}_1\) &
  \multicolumn{1}{r|}{104.3272} &
  \multicolumn{1}{r|}{101.1930} &
  \multicolumn{1}{r|}{107.4009} &
  \multicolumn{1}{r|}{100.1121} &
  \multicolumn{1}{r|}{\textbf{83.9749}} &
  99.4016 \\ \cline{2-8} 
 &
  \(\bm{K}_2\) &
  \multicolumn{1}{r|}{109.7507} &
  \multicolumn{1}{r|}{89.6992} &
  \multicolumn{1}{r|}{85.8899} &
  \multicolumn{1}{r|}{79.7213} &
  \multicolumn{1}{r|}{\textbf{77.6889}} &
  88.5500 \\ \hline
\end{tabular}%
}
\end{table}

Table \ref{tab: results} summarizes the results for both \(\bm{K}_1\) and \(\bm{K}_2\) such that time of flight is minimized. It can be observed that considering \(\bm{K}_2\) in the CLF improves the time of flight in each transfer case, with the improvement being more evident for maneuvers that occur over a longer time horizon. It can be observed that the global minimum is achievable for the shorter duration transfers, especially Cases A and C. This is evident by PSO being able to find roughly the same solution for all 5 runs of each weighting matrix. It is evident that the global solution is effected by the type of parameterization for these cases. 

Multiple local minima are found for the other transfer cases. In particular, Case E  demonstrates the existence of many local minima, with the times of flight ranging broadly and with the worst solution achieved being found by the full weighting matrix, whereas in the other cases the worst solution was always found by the diagonal weighting matrix. On the other hand, the best solution is found by the full matrix representation in each transfer case and the average times of flights for the full weighting matrix \(\bm{K}_2\) are lower than those from the diagonal weighting matrix \(\bm{K}_1\) in each transfer case. 

\begin{figure}
    \centering
    \includegraphics[width=0.8\linewidth]{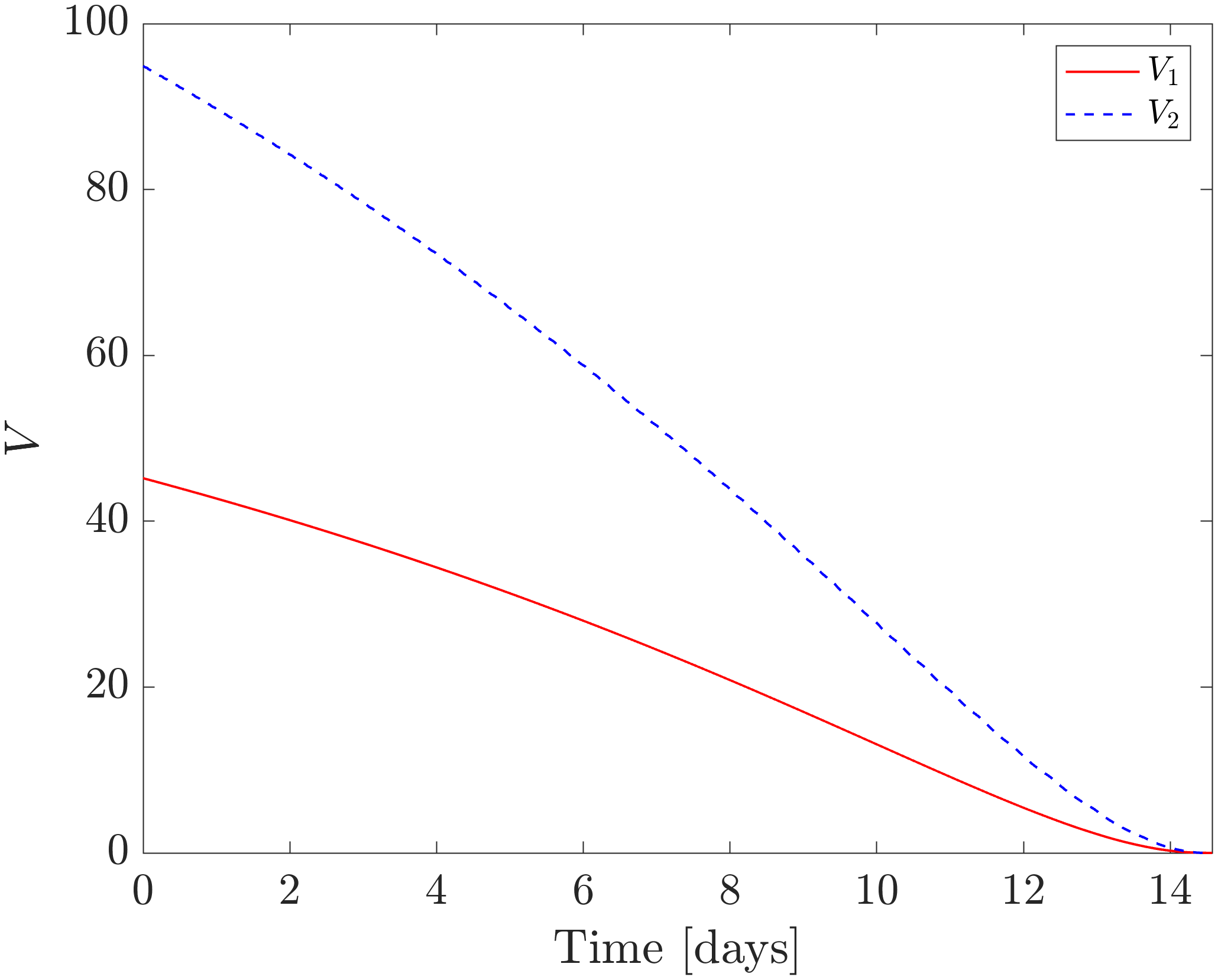}
    \includegraphics[width=0.8\linewidth]{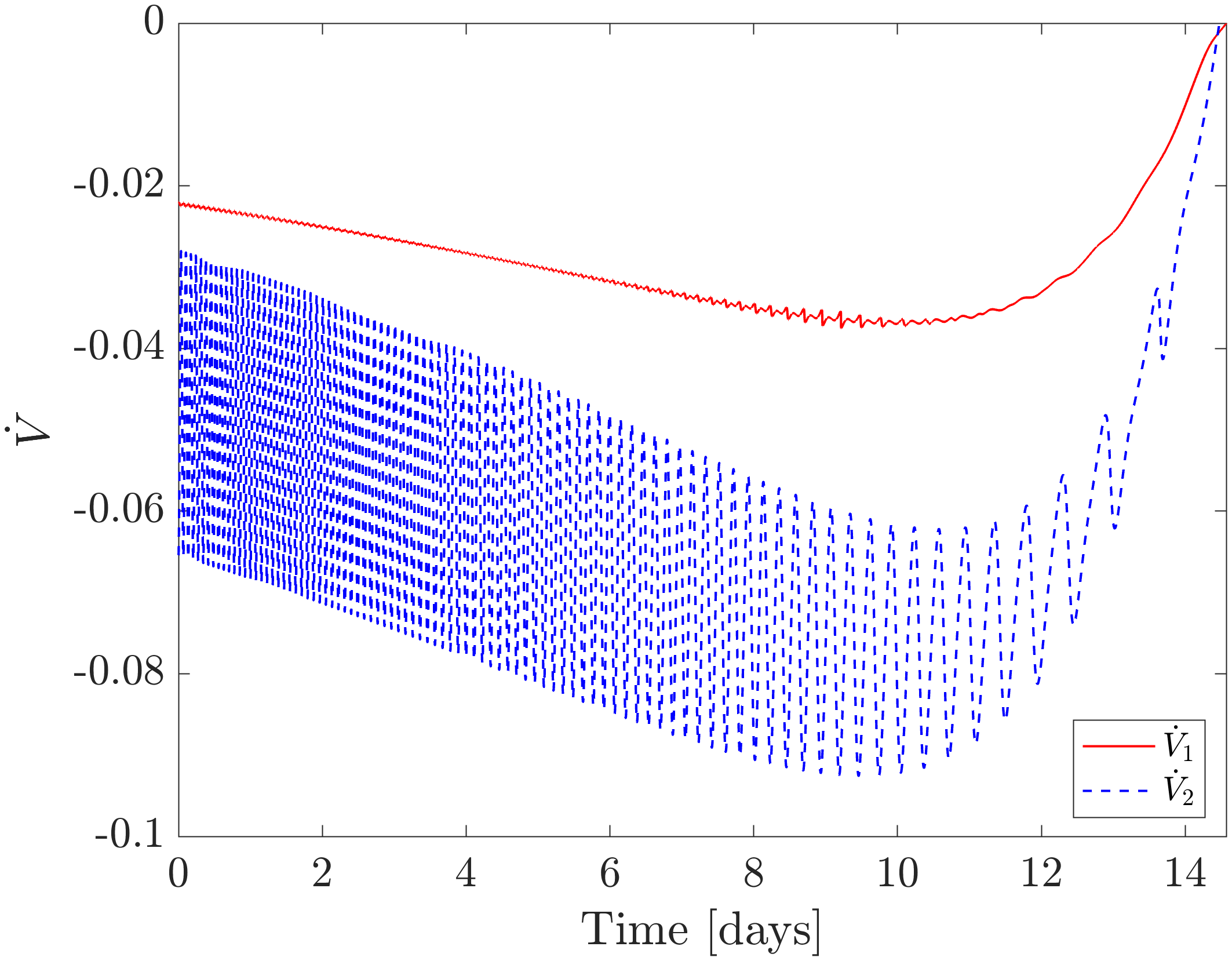}
    \caption{Case A: Lyapunov function and its time-derivative vs. time.}
    \label{fig: case a V}
\end{figure}

\begin{figure}
    \centering
    \includegraphics[width=0.8\linewidth]{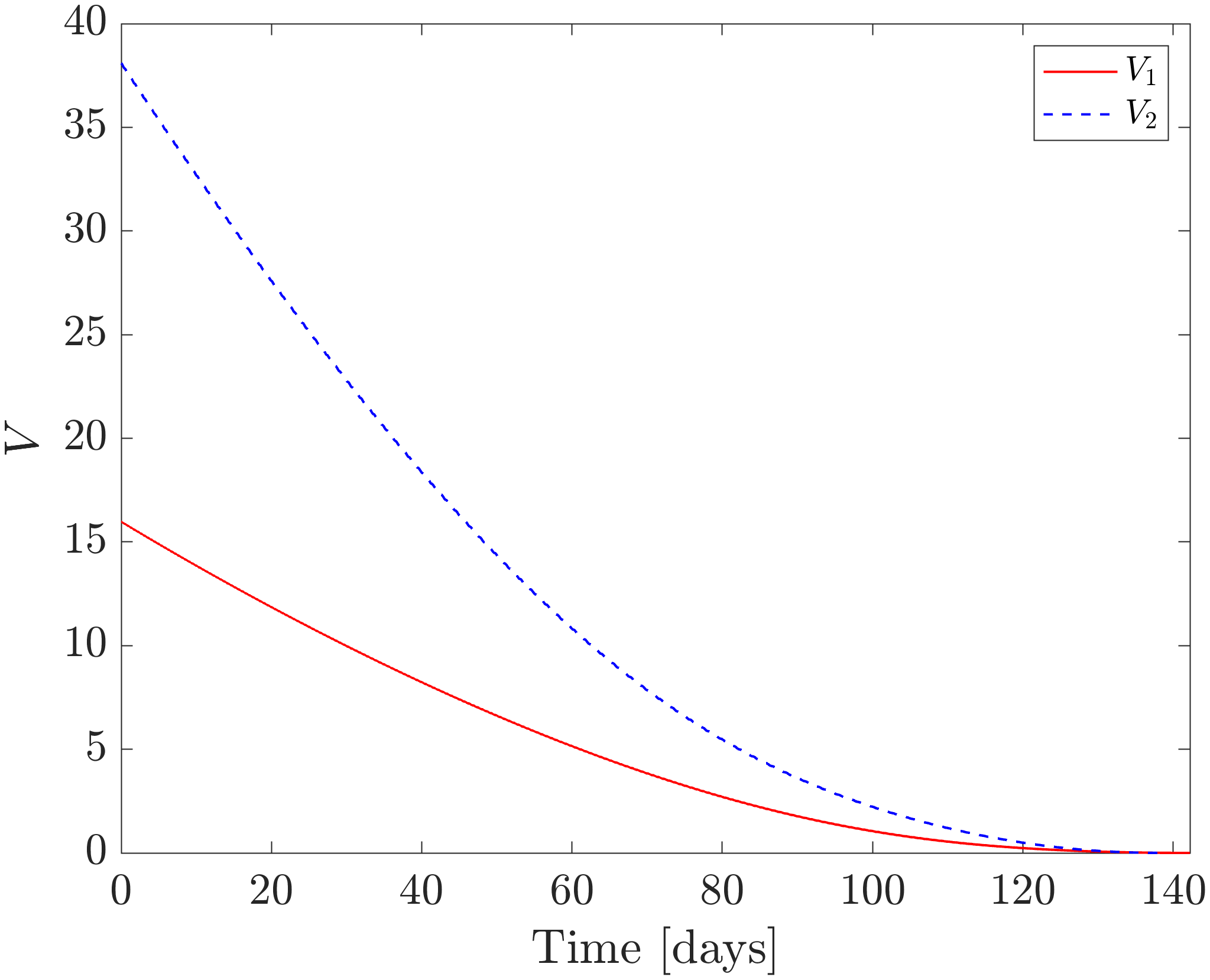}
    \includegraphics[width=0.8\linewidth]{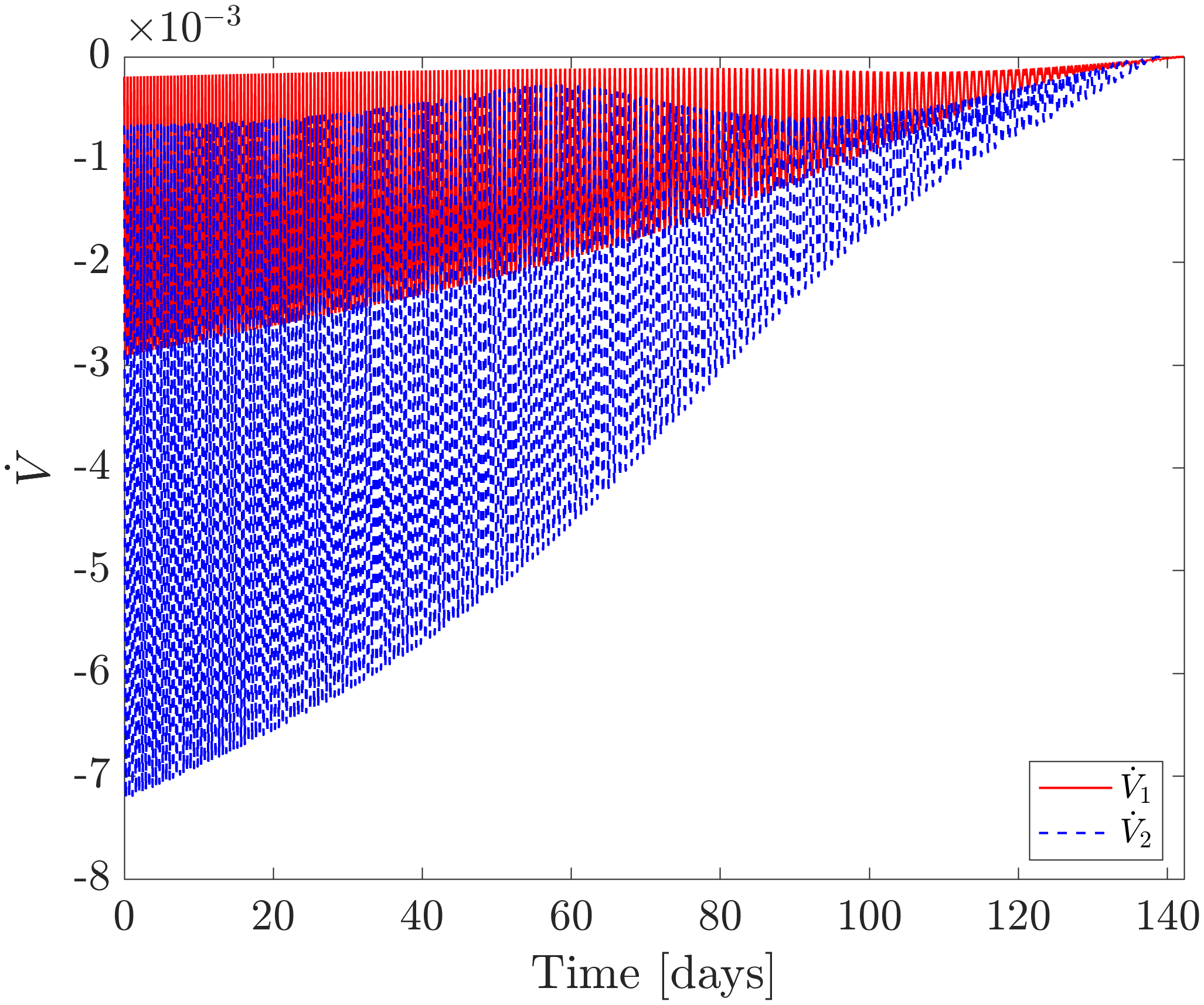}
    \caption{Case B: Lyapunov function and its time-derivative vs. time.}
    \label{fig: case b V}
\end{figure}

\begin{figure}
    \centering
    \includegraphics[width=0.8\linewidth]{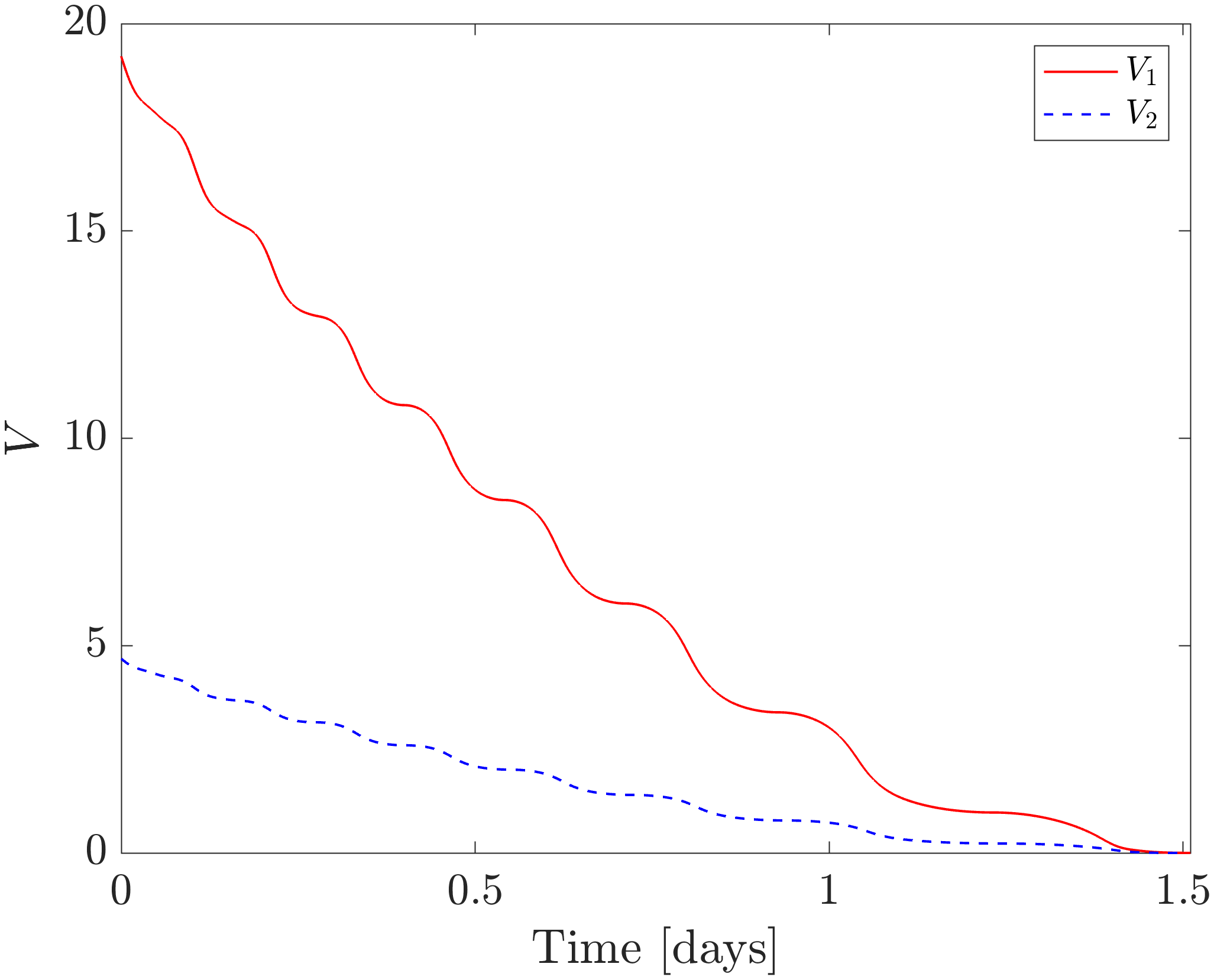}
    \includegraphics[width=0.8\linewidth]{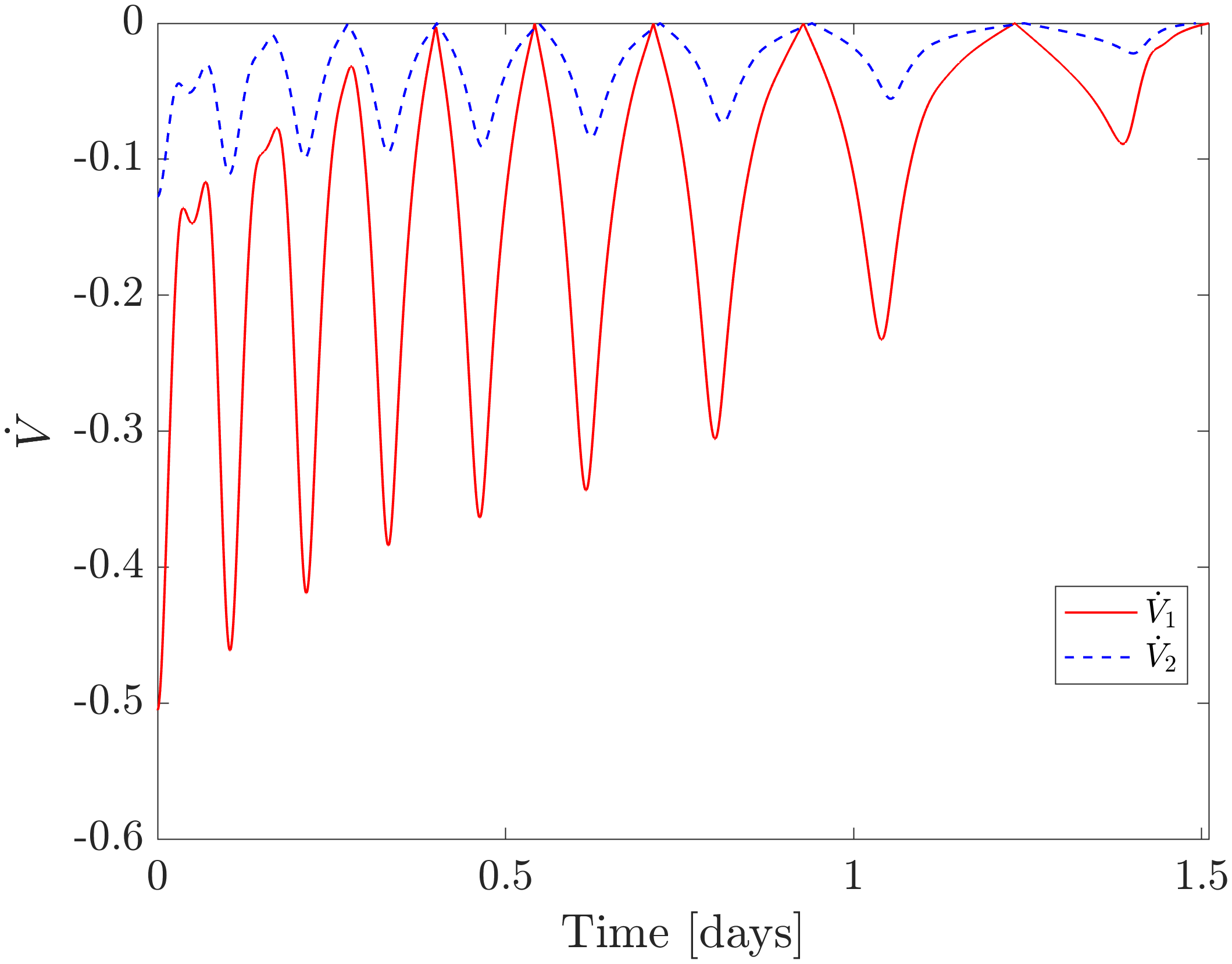}
    \caption{Case C: Lyapunov function and its time-derivative vs. time.}
    \label{fig: case c V}
\end{figure}

\begin{figure}
    \centering \hspace{7mm}
    \includegraphics[width=0.8\linewidth]{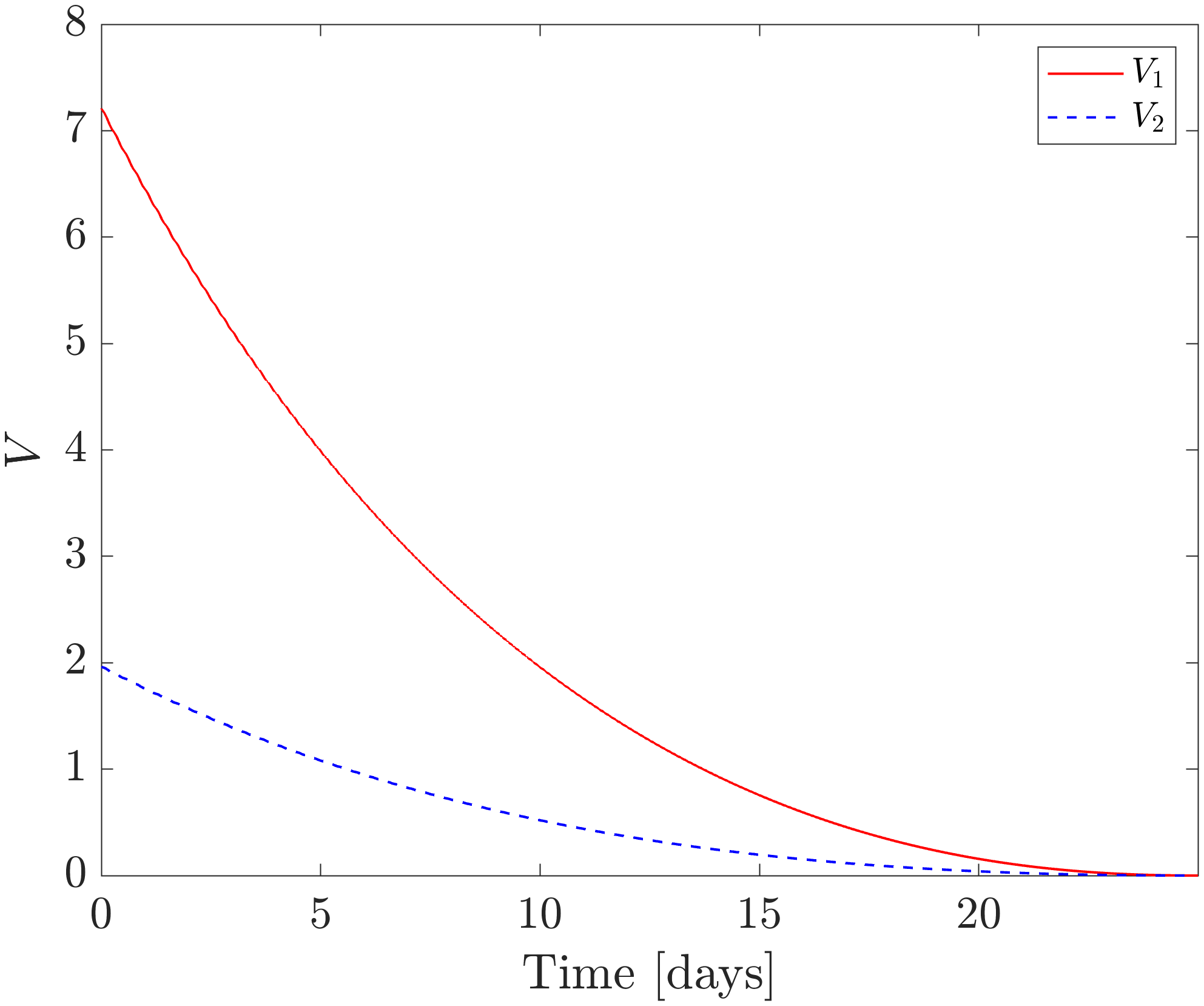}
    \includegraphics[width=0.88\linewidth]{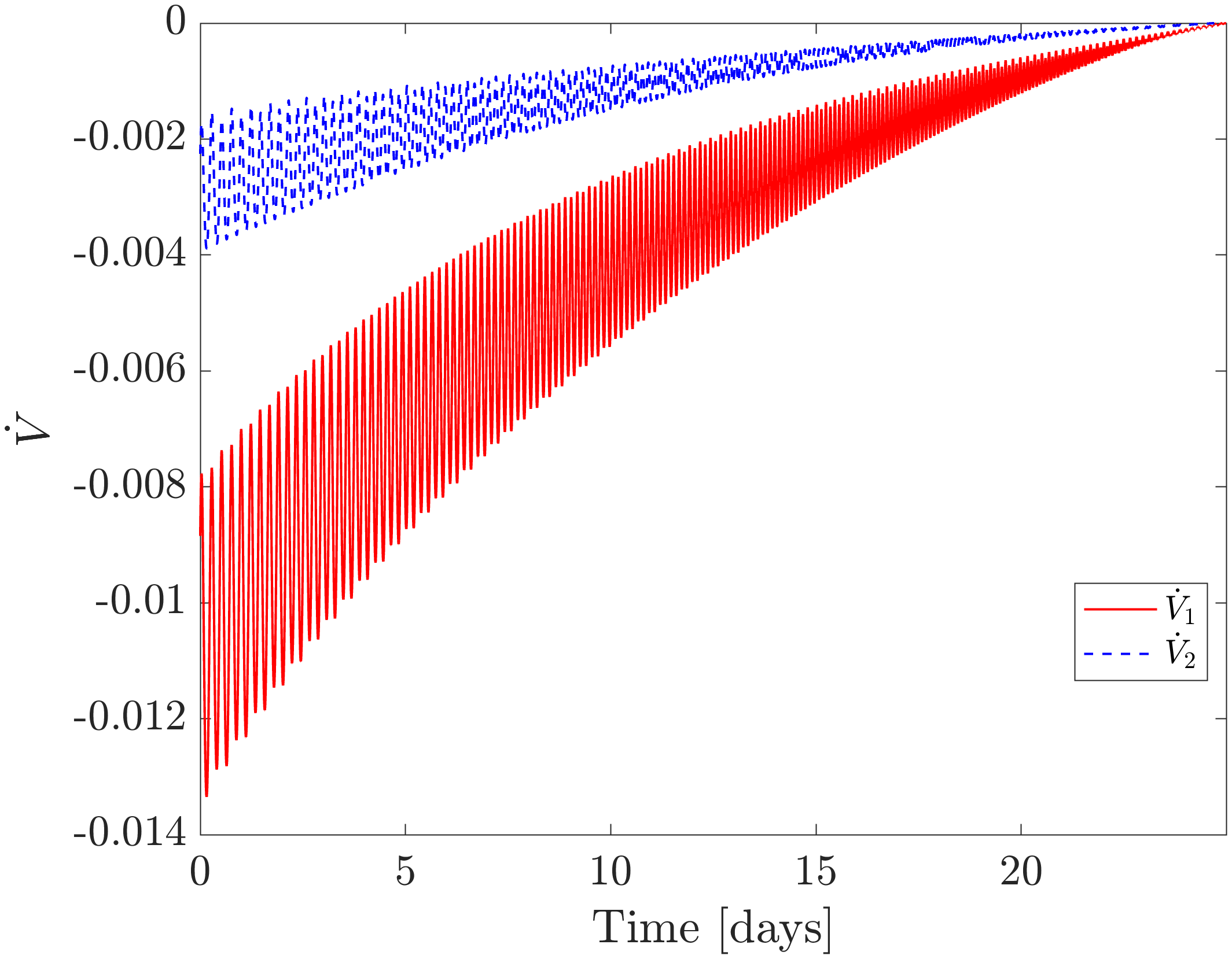}
    \caption{Case D: Lyapunov function and its time-derivative vs. time.}
    \label{fig: case d V}
\end{figure}

\begin{figure}
    \centering
    \includegraphics[width=0.8\linewidth]{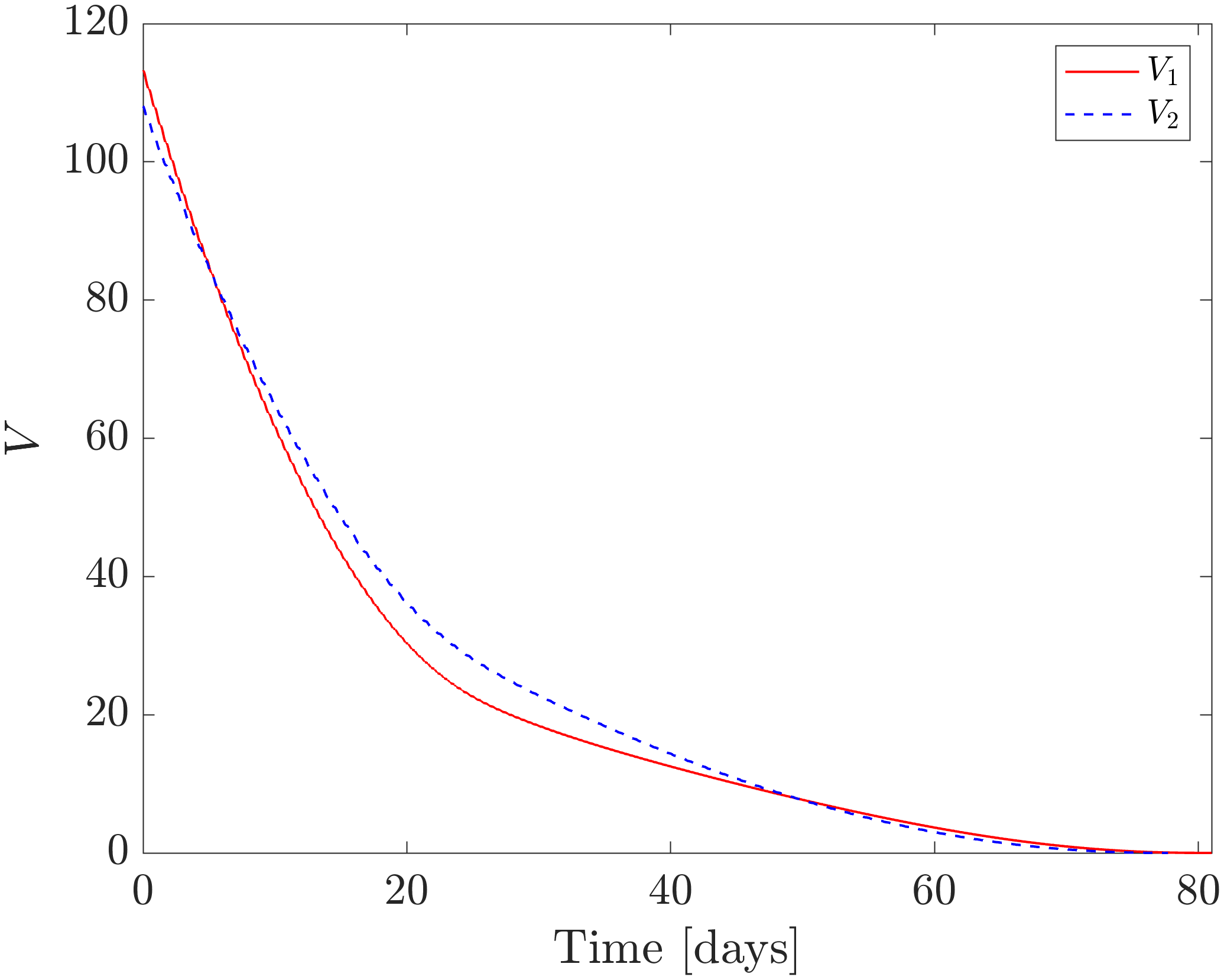}
    \includegraphics[width=0.8\linewidth]{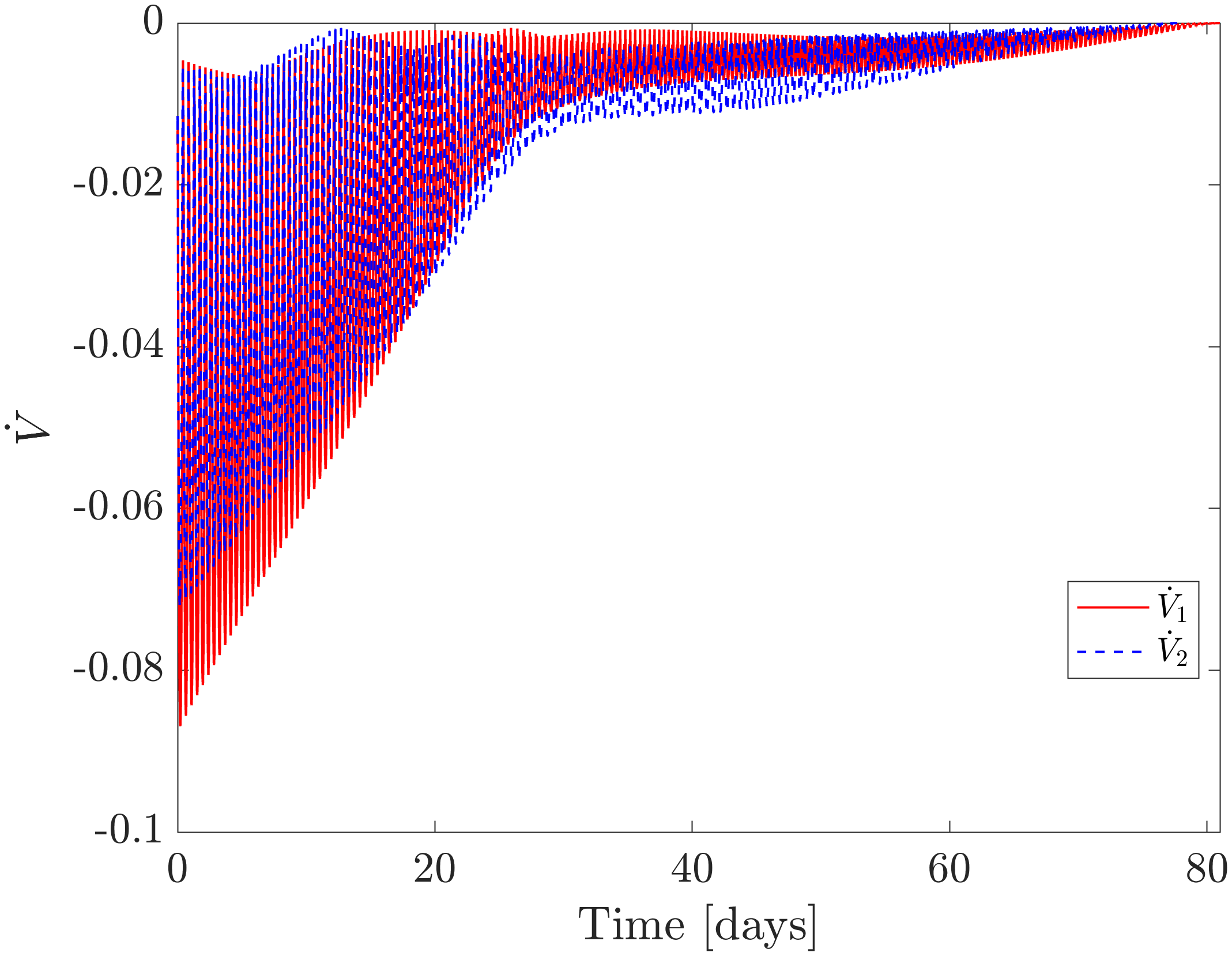}
    \caption{Case E: Lyapunov function and its time-derivative vs. time.}
    \label{fig: case e V}
\end{figure}

Note that in all solution plots, subscript `1' corresponds to the solutions from the diagonal weighting matrix, \(\bm{K}_1\), and subscript `2' corresponds to solutions from the full weighting matrix, \(\bm{K}_2\). Plots of the Lyapunov function and its time derivative versus time, for the best solution from each weighting matrix, are shown in Figures \ref{fig: case a V}, \ref{fig: case b V}, \ref{fig: case c V}, \ref{fig: case d V}, and \ref{fig: case e V}. It can be observed that the CLF in fact remains positive-definite in every transfer case. However, it is not intuitively obvious in these plots how the full weighting matrix \(\bm{K}_2\) is driving the CLF to 0 faster than the diagonal weighting matrix \(\bm{K}_1\). 

\begin{figure}
    \centering
    \includegraphics[width=0.8\linewidth]{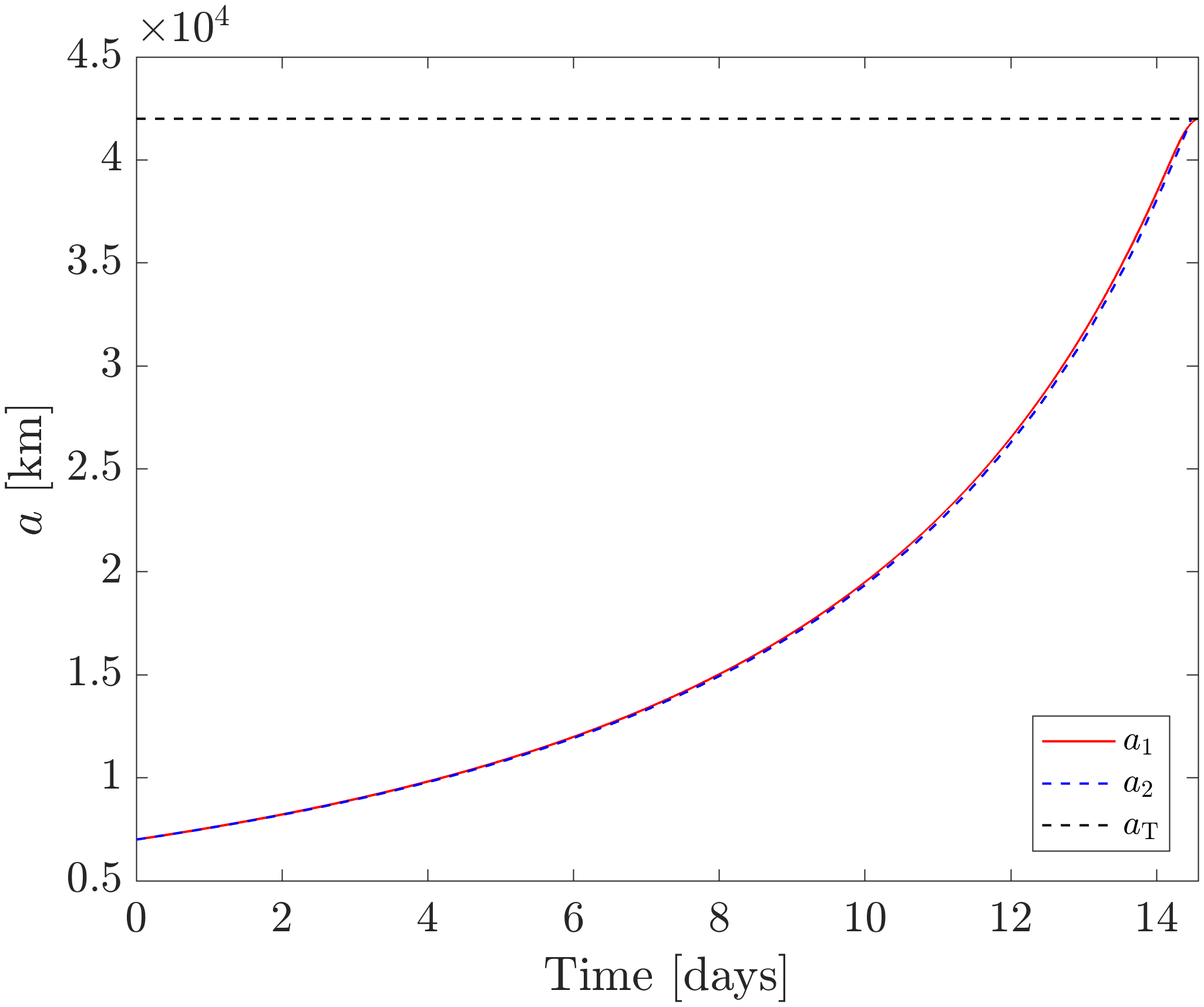}
    \includegraphics[width=0.8\linewidth]{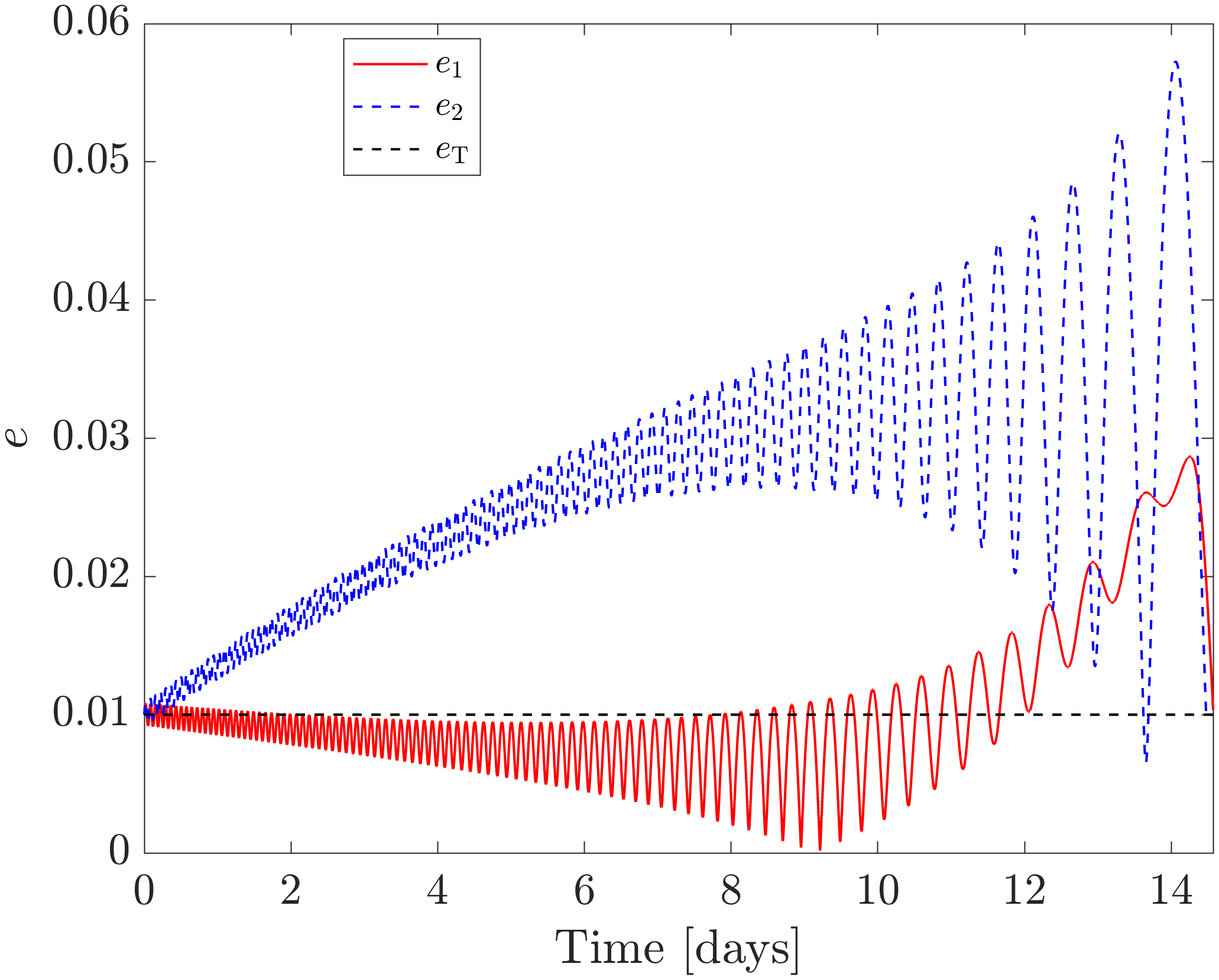}
    \caption{Case A: Classical orbital elements ($a$ and $e$) vs. time.}
    \label{fig: case a coe}
\end{figure}

\begin{figure}
    \centering
    \includegraphics[width=0.49\linewidth]{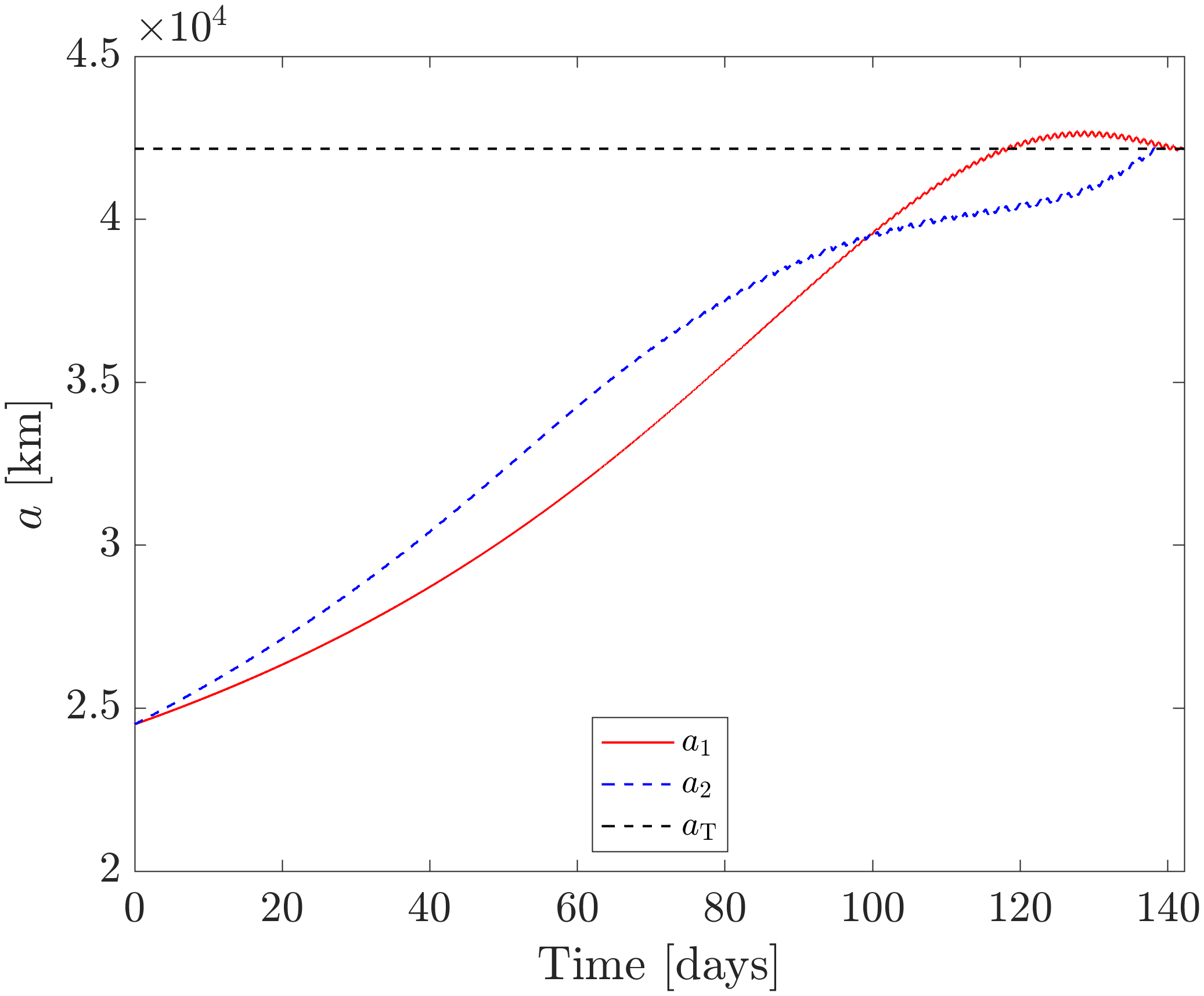}
    \includegraphics[width=0.49\linewidth]{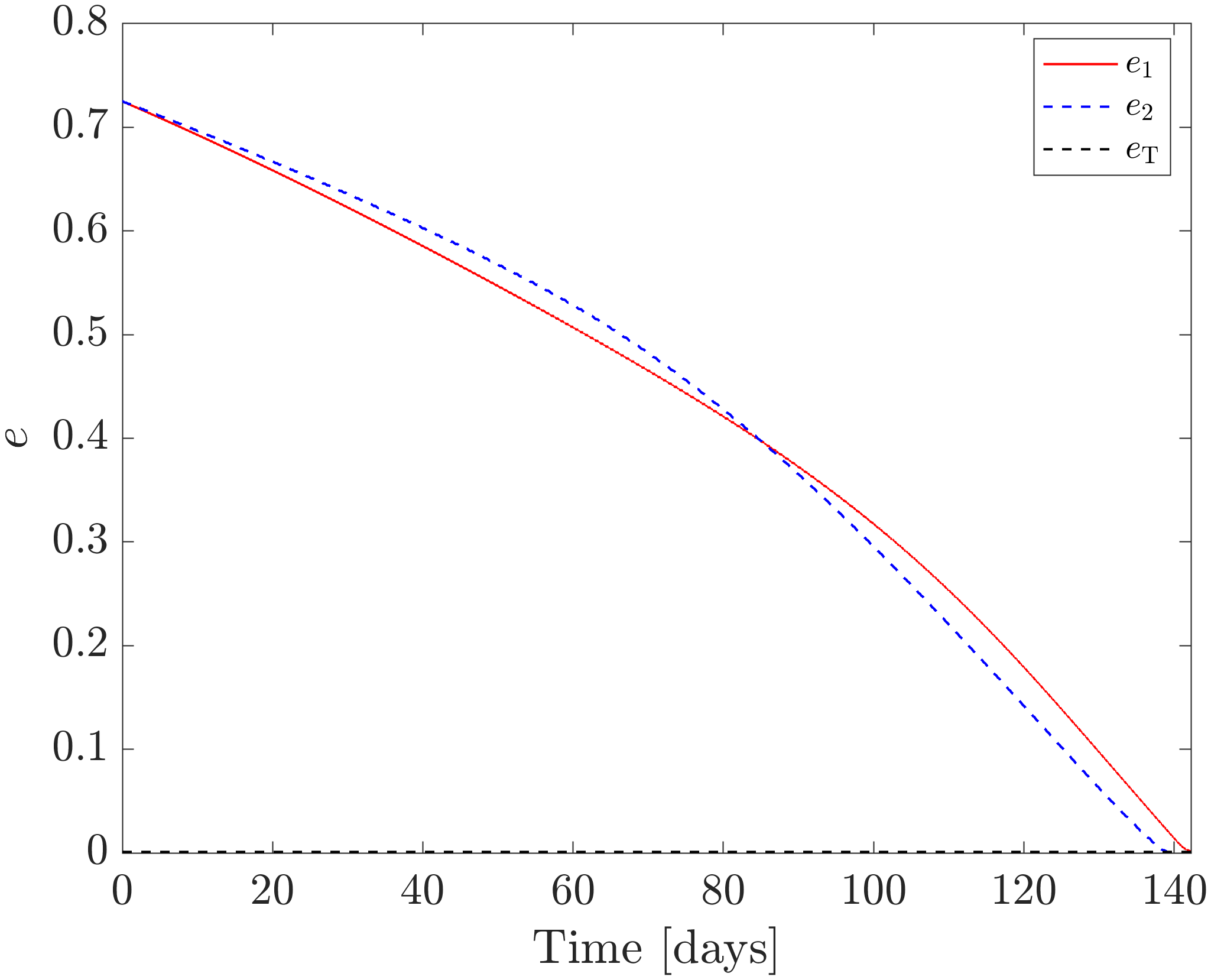}
    \includegraphics[width=0.49\linewidth]{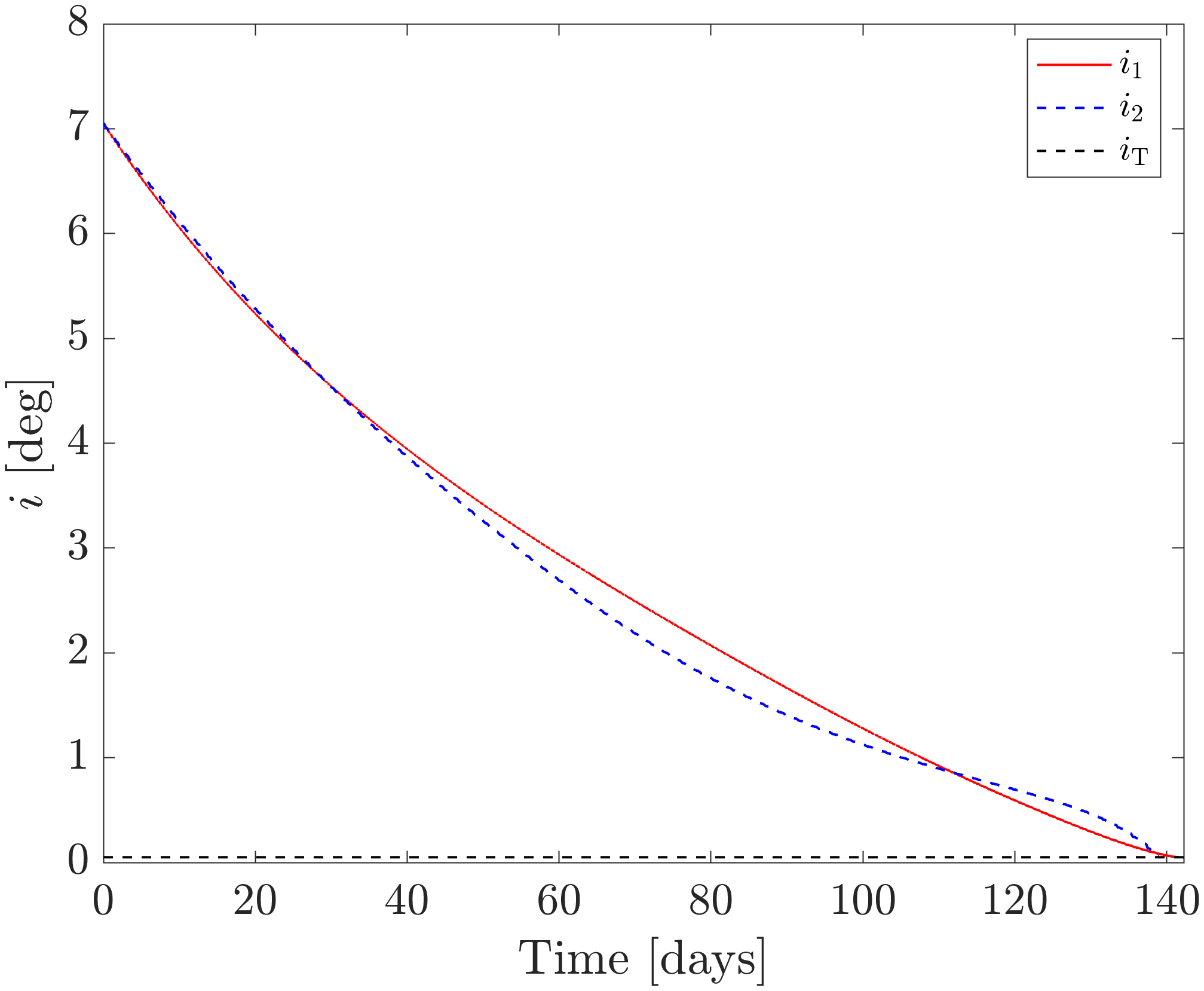}
    \caption{Case B: Classical orbital elements ($a$, $e$, and $i$) vs. time.}
    \label{fig: case b coe}
\end{figure}

\begin{figure}
    \centering
    \includegraphics[width=0.49\linewidth]{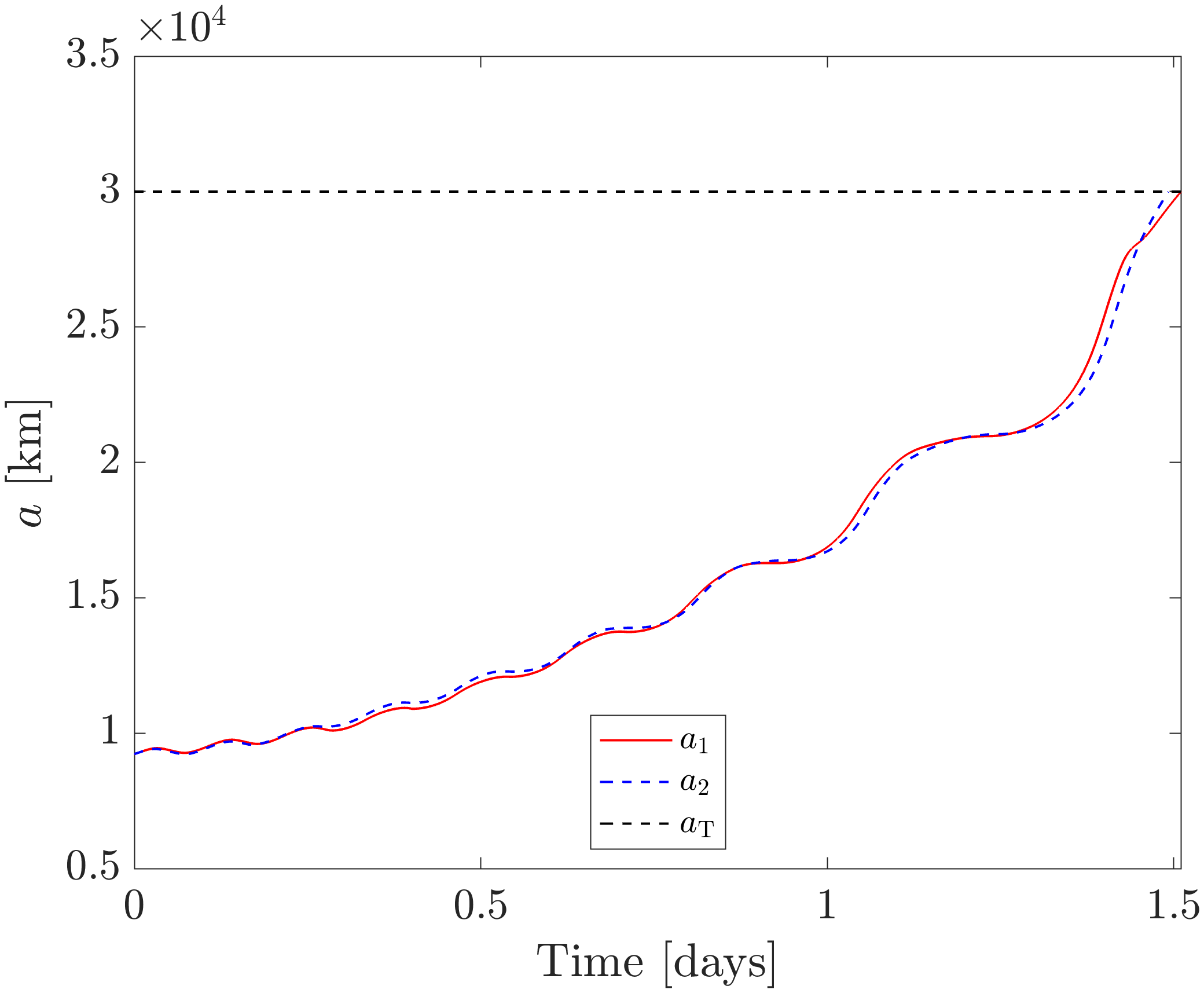}
    \includegraphics[width=0.49\linewidth]{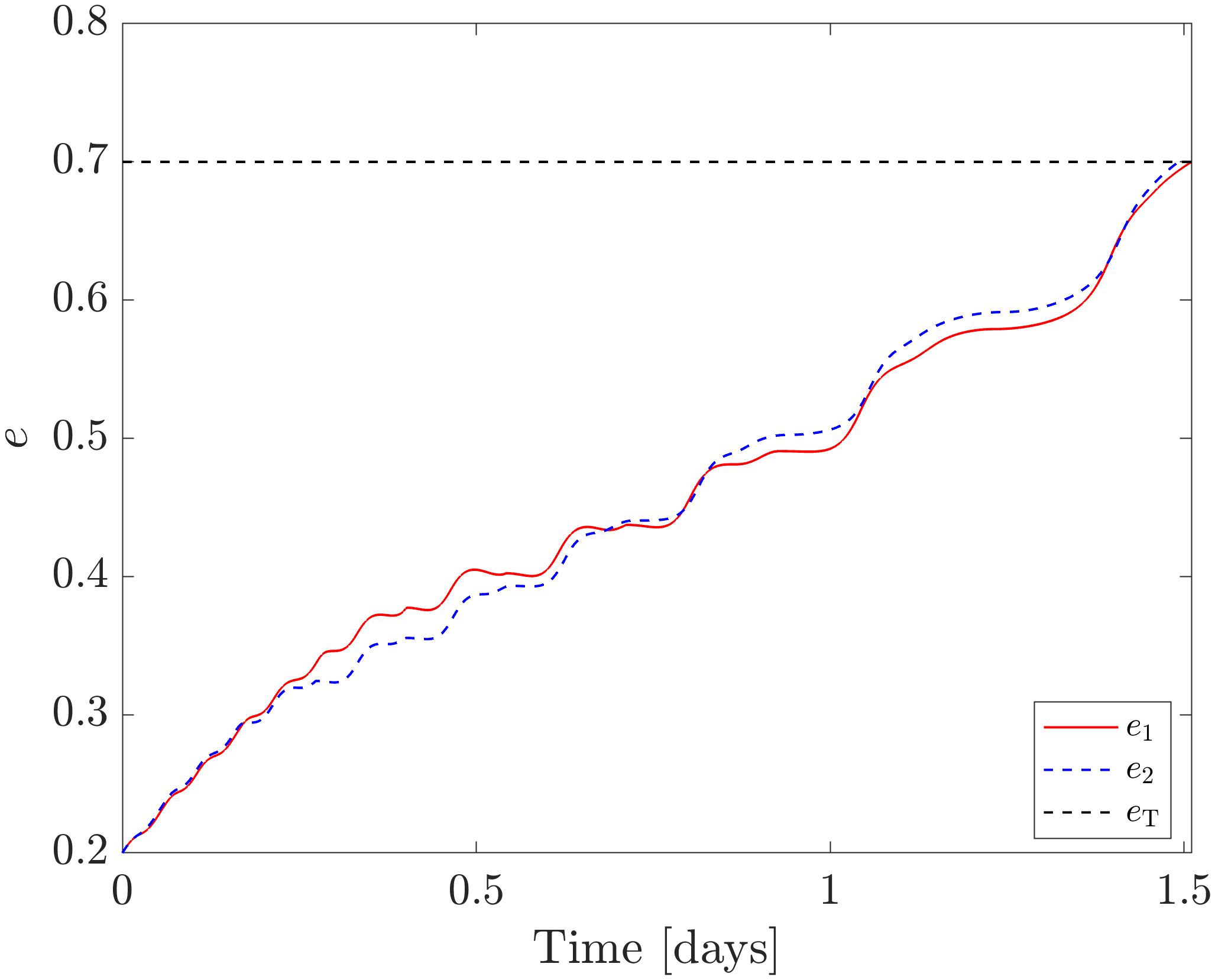}
    \caption{Case C: Classical orbital elements ($a$ and $e$) vs. time.}
    \label{fig: case c coe}
\end{figure}

\begin{figure}
    \centering
    \includegraphics[width=0.49\linewidth]{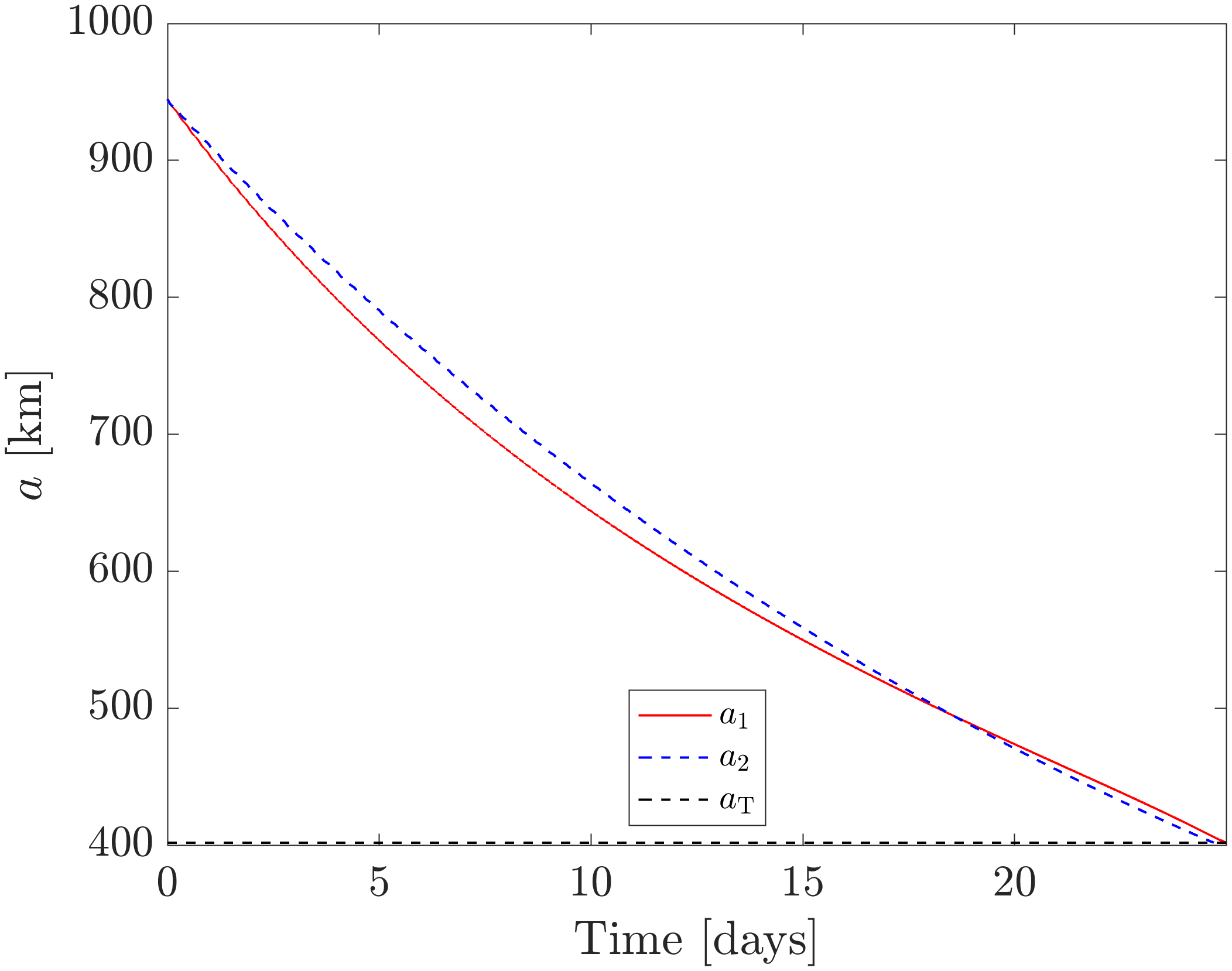}
    \includegraphics[width=0.49\linewidth]{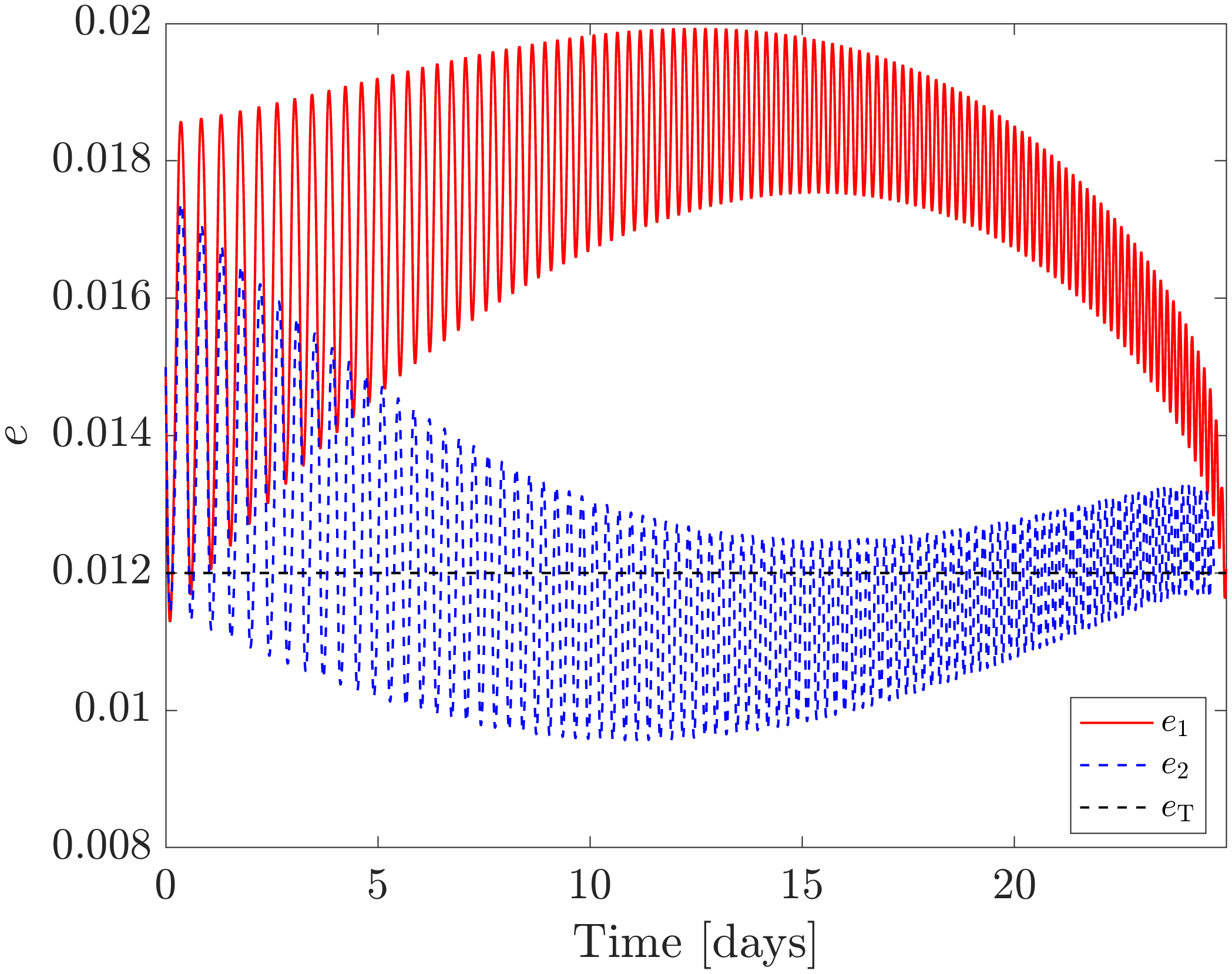}
    \includegraphics[width=0.49\linewidth]{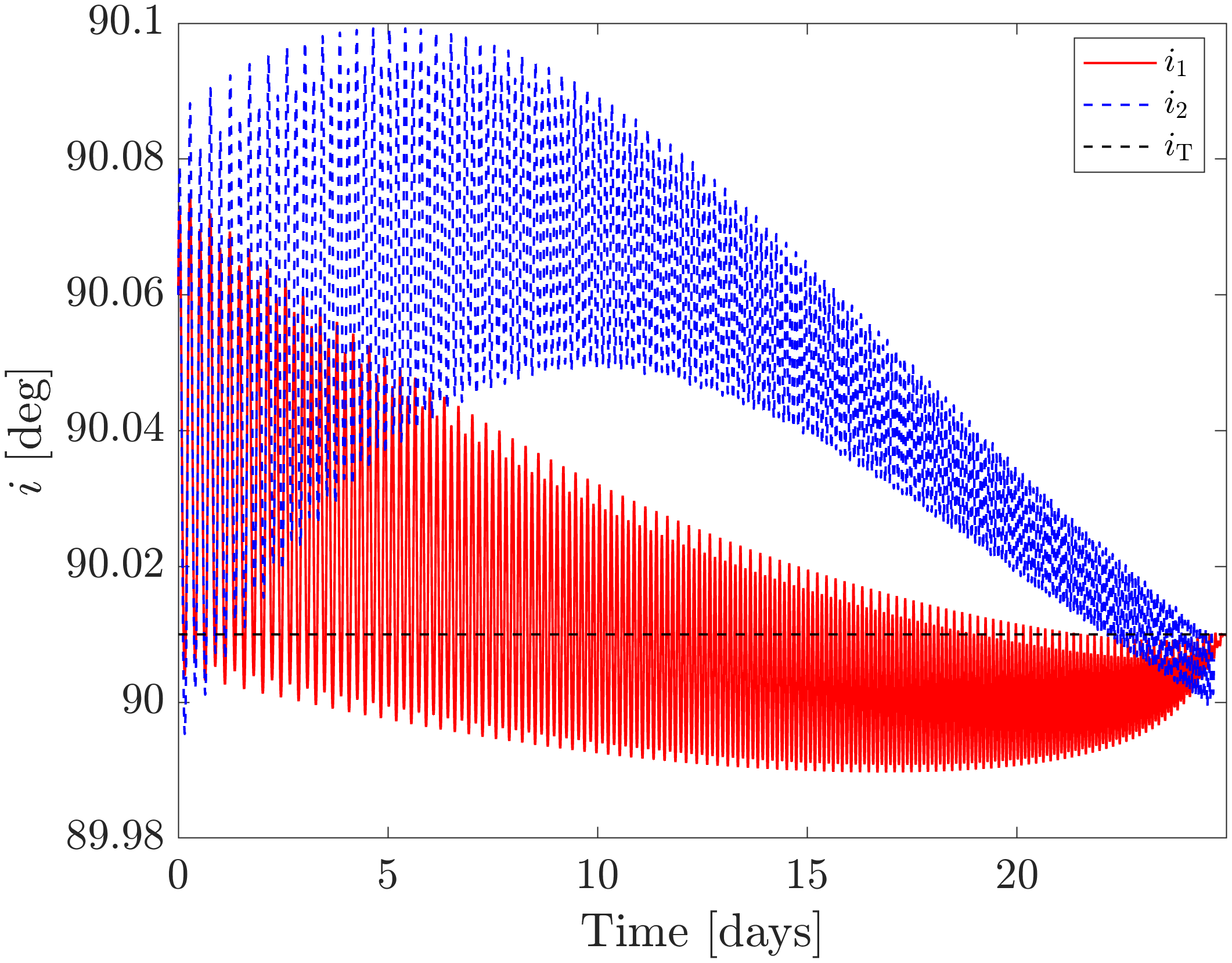}
    \includegraphics[width=0.49\linewidth]{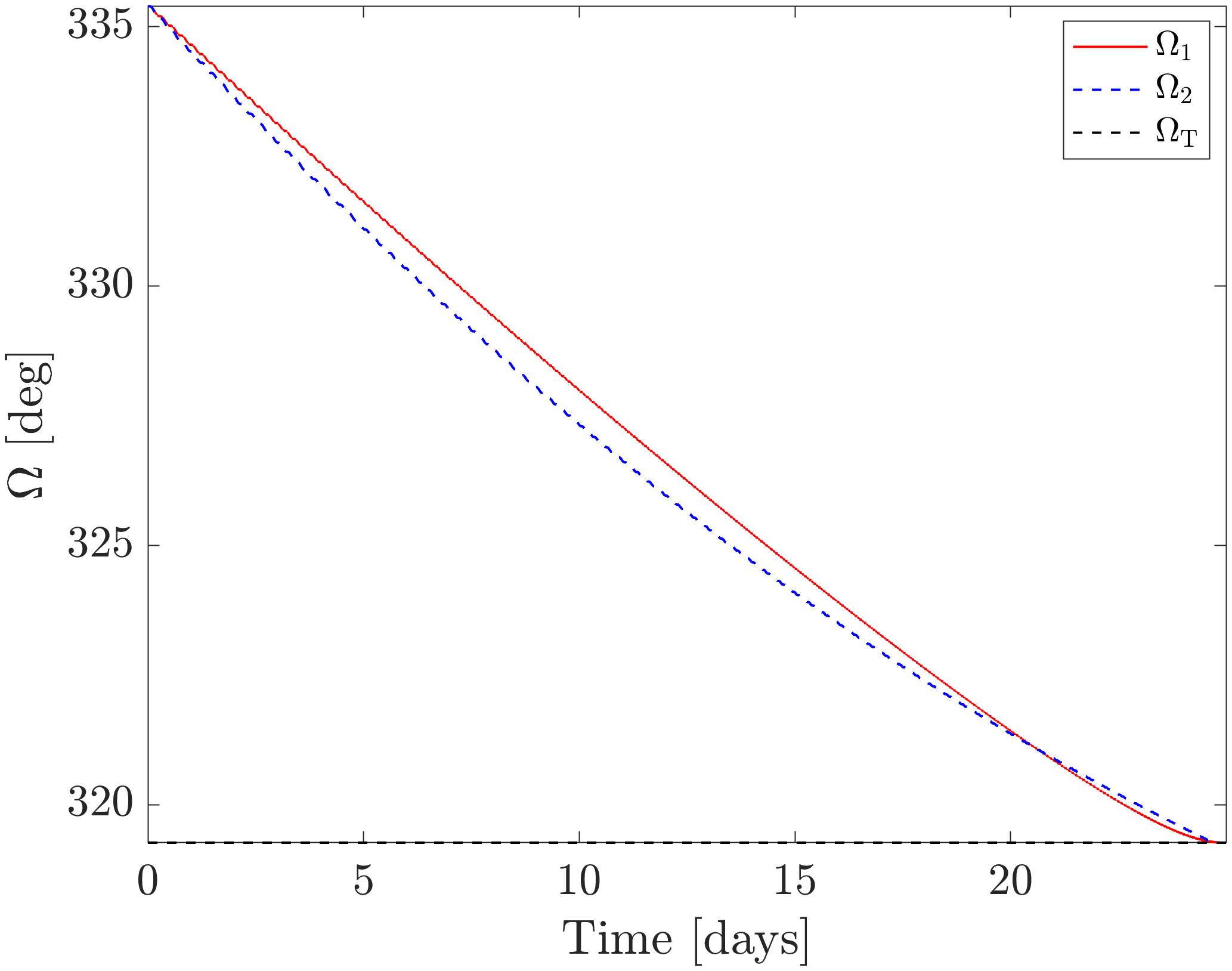}
    \caption{Case D: Classical orbital elements ($a$, $e$, $i$, and $\Omega$) vs. time.}
    \label{fig: case d coe}
\end{figure}

\begin{figure}
    \centering
    \includegraphics[width=0.49\linewidth]{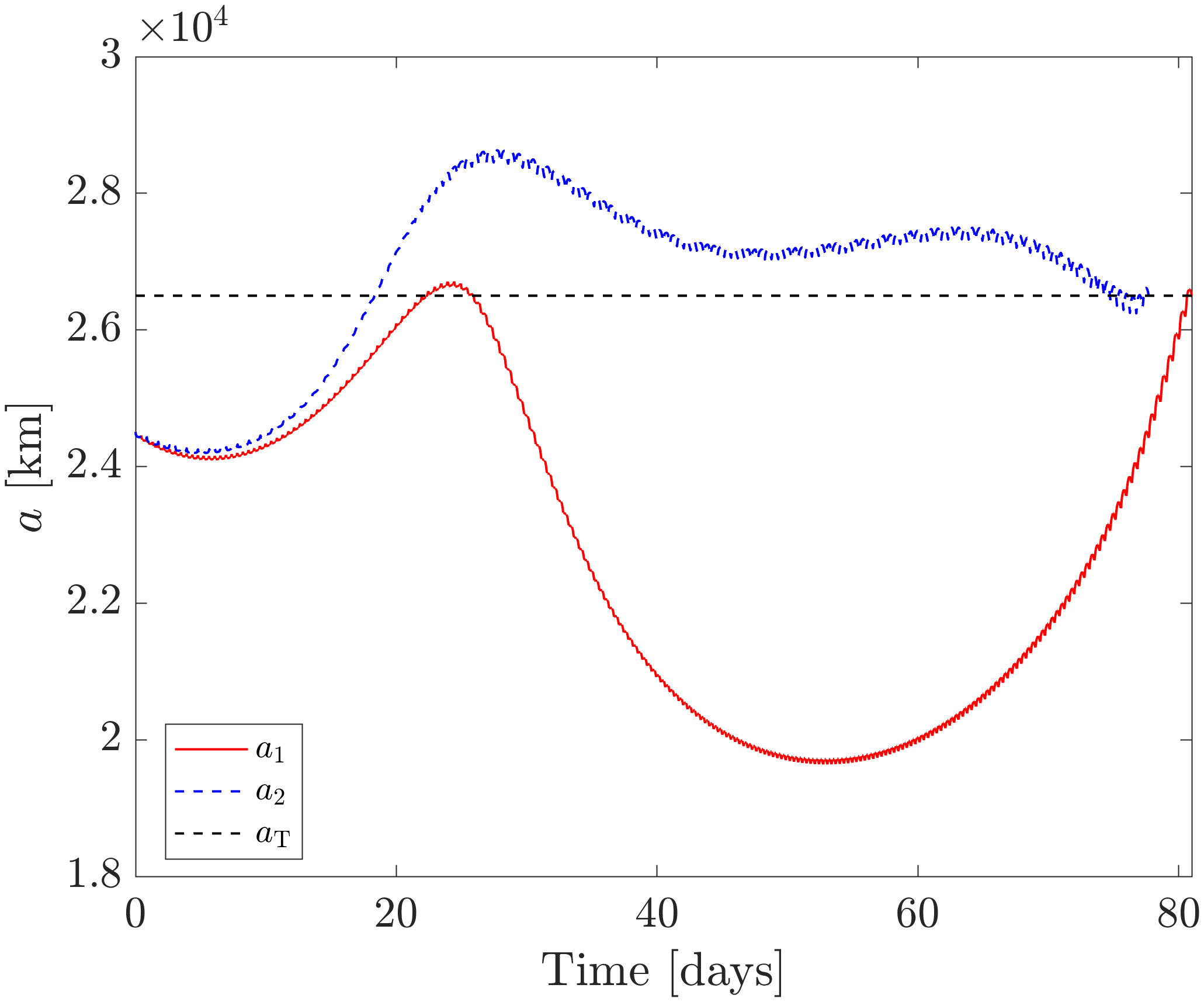}
    \includegraphics[width=0.49\linewidth]{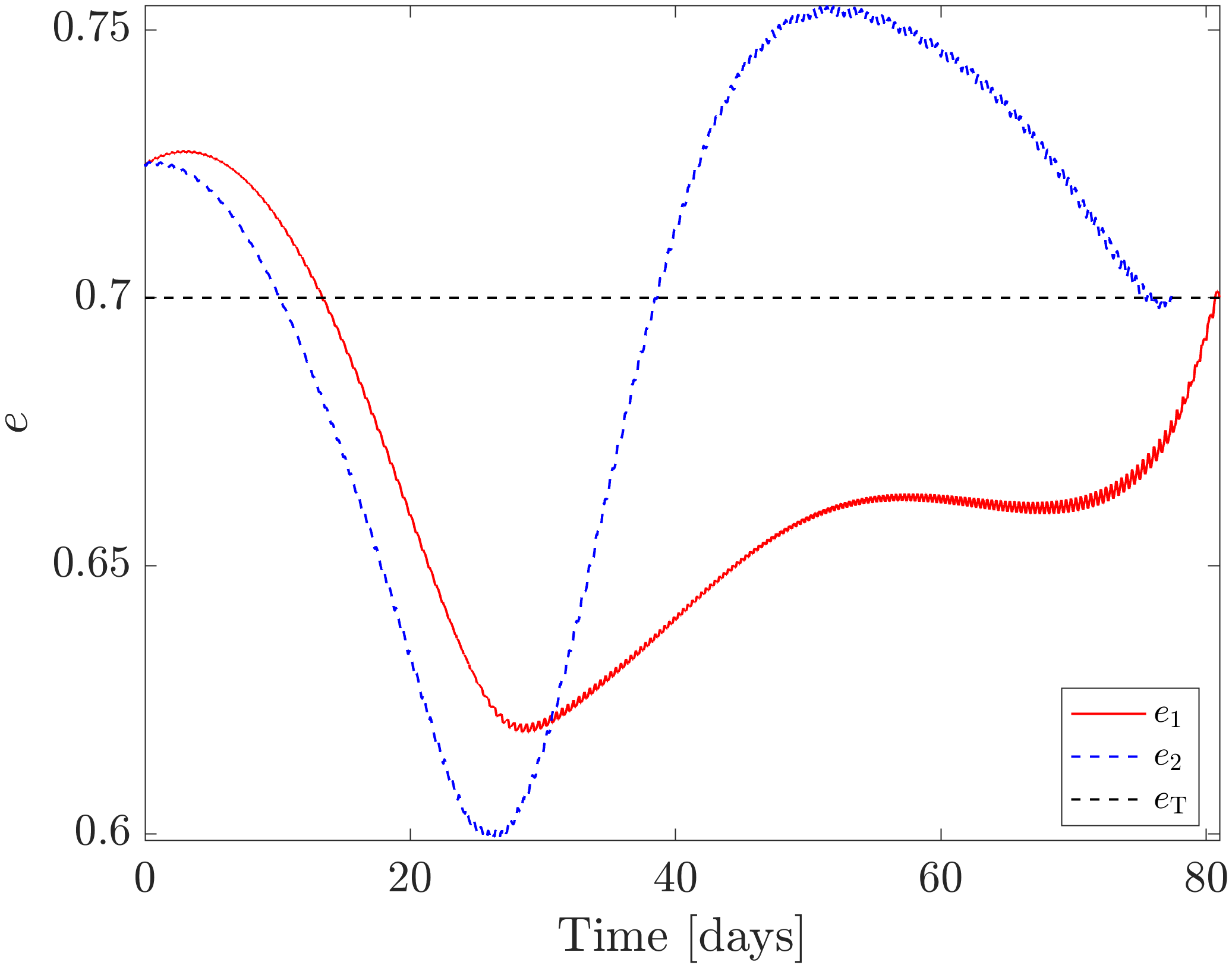}
    \includegraphics[width=0.49\linewidth]{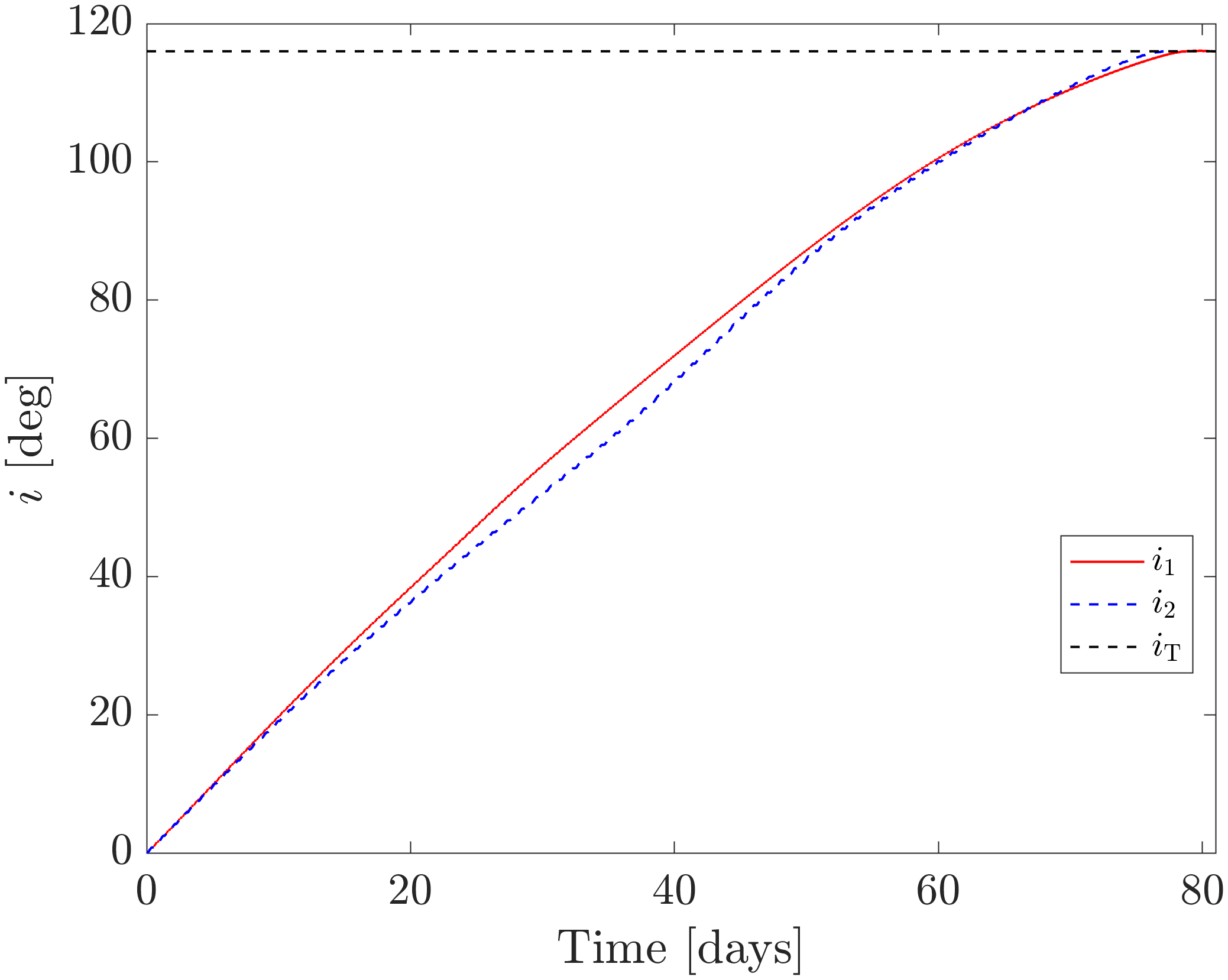}
    \includegraphics[width=0.49\linewidth]{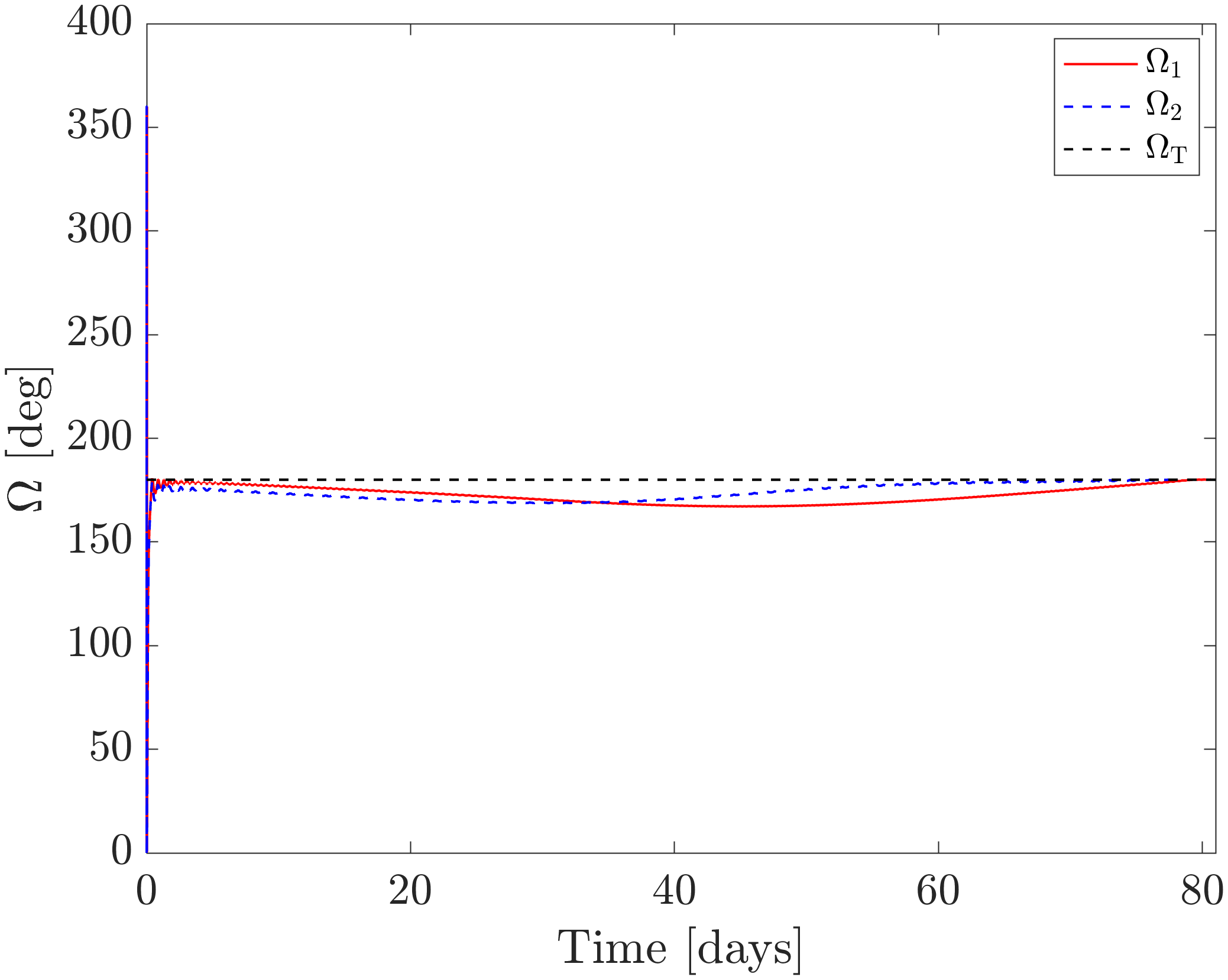}
    \includegraphics[width=0.49\linewidth]{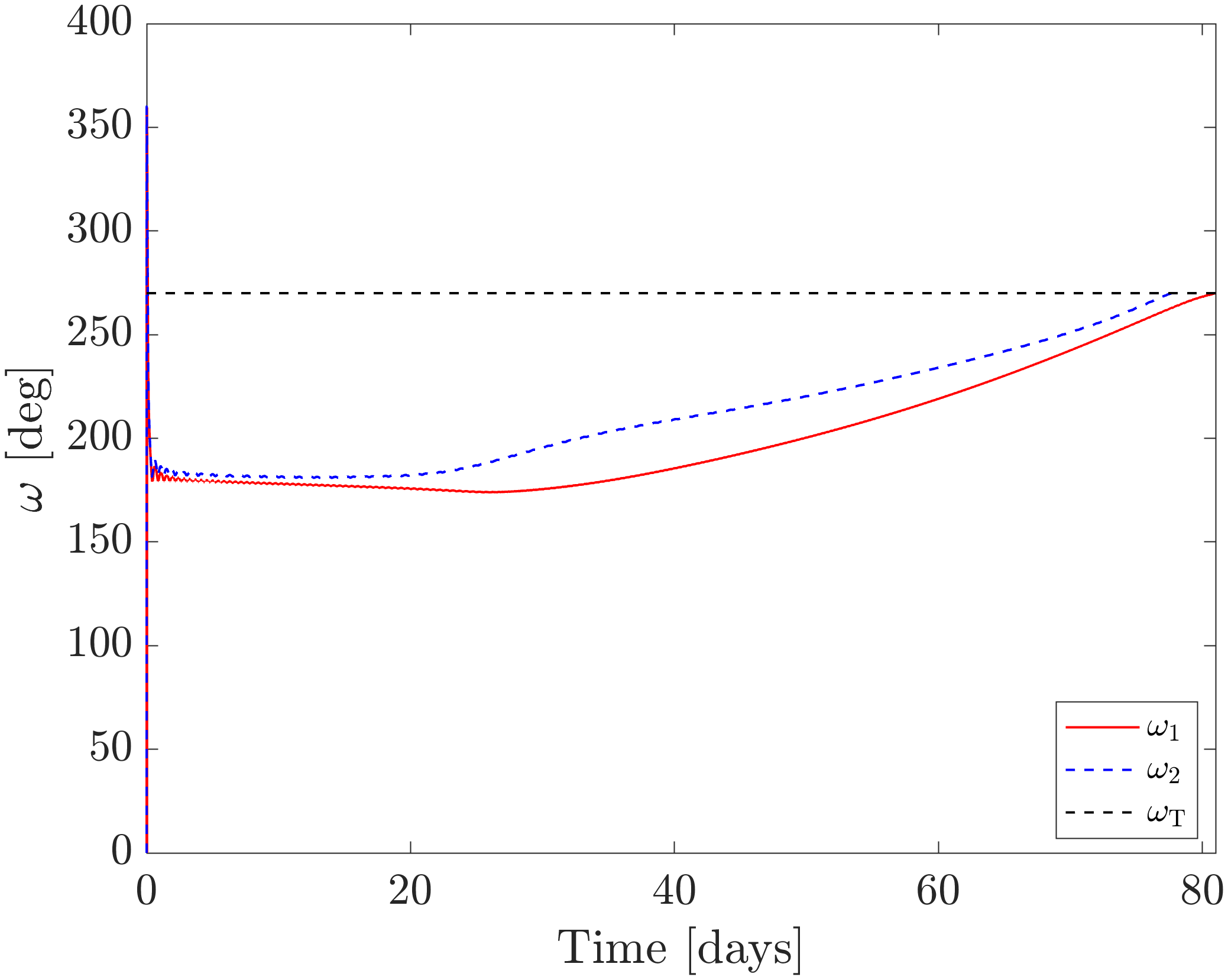}
    \caption{Case E: Classical orbital elements ($a$, $e$, $i$, $\Omega$, and $\omega$) vs. time.}
    \label{fig: case e coe}
\end{figure}

Plots of the time histories of the classical orbital elements for the best solution from each weighting matrix are shown in Figures \ref{fig: case a coe}, \ref{fig: case b coe}, \ref{fig: case c coe}, \ref{fig: case d coe}, and \ref{fig: case e coe}. It can be observed that the full weighting matrix, \(\bm{K}_2\), drives the orbital elements to their target values differently (and likely more efficiently based on the improved times of flight) than the diagonal weighting matrix \(\bm{K}_1\). The orbital elements are driven to its target value very differently in some cases like, for instance, the semi-major axis and eccentricity for Case E in Figure \ref{fig: case e coe}. The \(\bm{K}_2\) parameterization, taking into account the cross-coupling errors of the orbital elements, is likely what results in a solution with a better time of flight than the diagonal matrix \(\bm{K}_1\). 

\begin{figure}
    \centering
    \includegraphics[width=0.65\linewidth]{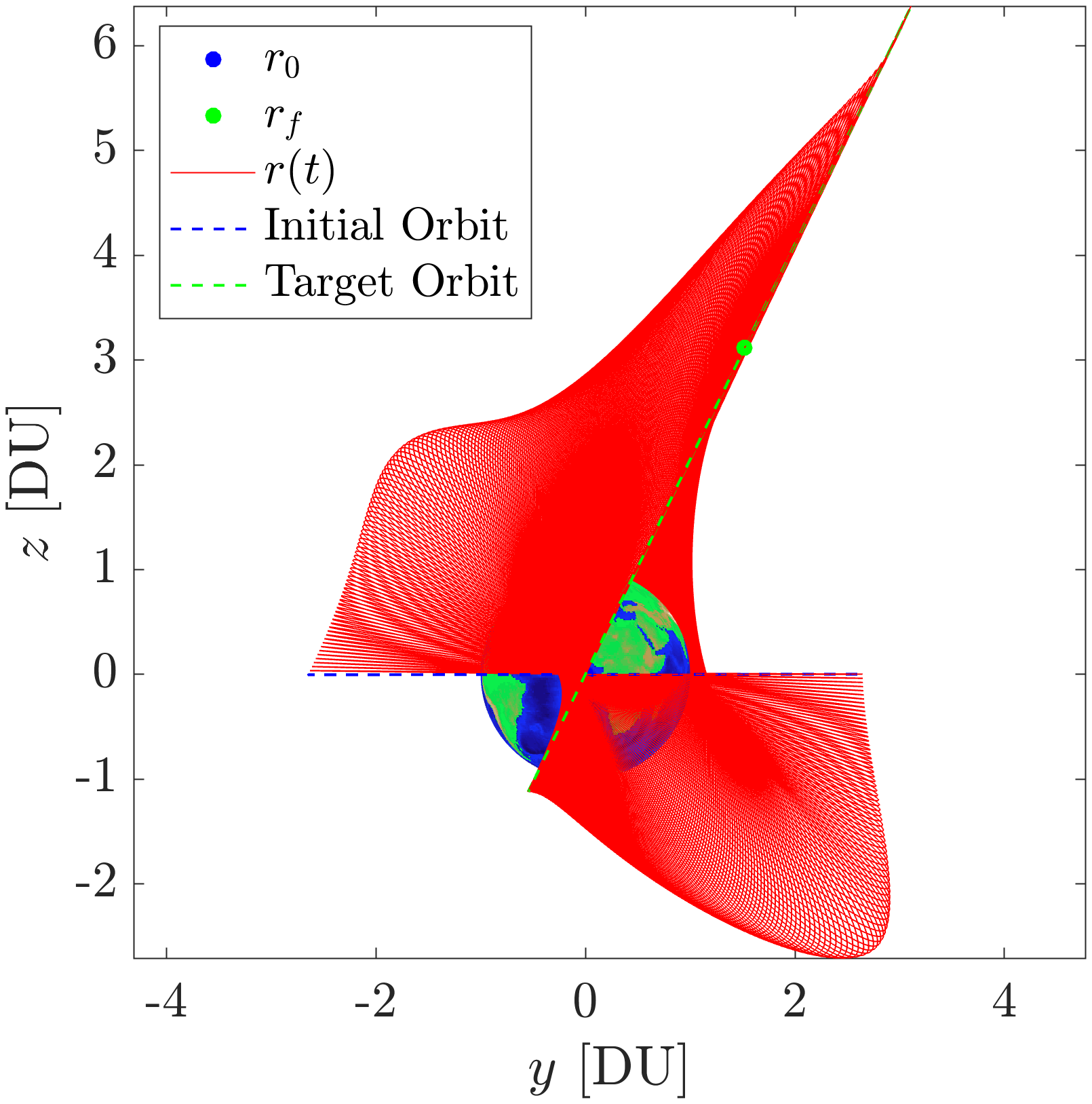}
    \includegraphics[width=0.65\linewidth]{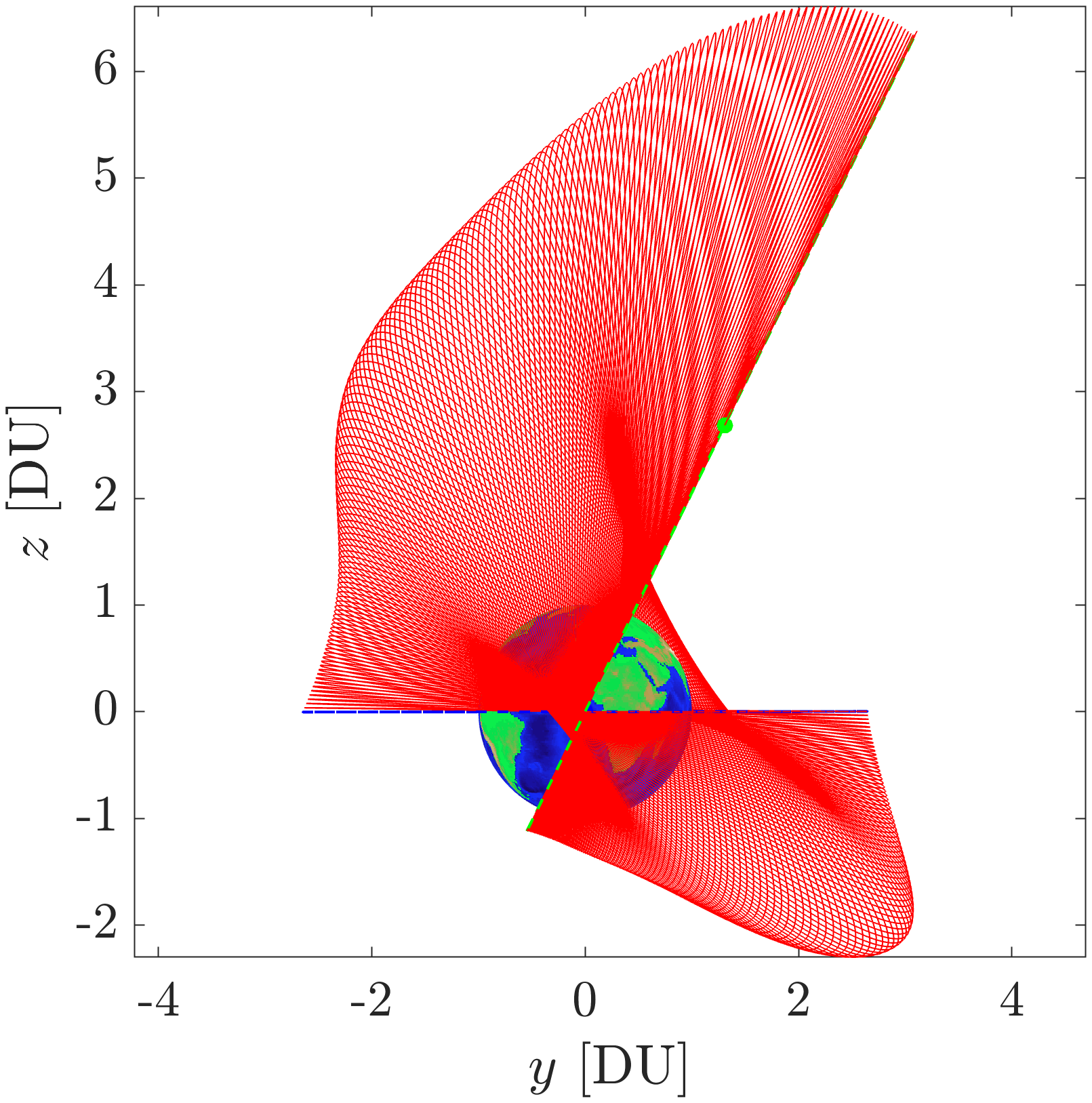}
    \caption{Case E: Trajectories for the best solutions obtained using the diagonal weighting matrix \(\bm{K}_1\) (top) and the full weighting matrix \(\bm{K}_2\) (bottom).}
    \label{fig: case e traj}
\end{figure}

Most of the solutions with the two weighting matrices result in trajectories that are not visually different. However, Figure \ref{fig: case e traj} shows a comparison of the best trajectories from each weighting matrix for Case E, which shows a significant difference. The inclination change follows a similar trend in both solutions. For \(\bm{K}_1\), the semi-major axis has first increased to beyond the semi-major axis of the target orbit, but then, it is reduced noticeably below the target value. On the other hand, the eccentricity is first decreased and then gradually increased. For \(\bm{K}_2\), the semi major axis is increased and maintained more or less until the end of the maneuver. The eccentricity, however, has changed in a sinusoidal manner. 

The weighting matrices corresponding to these two runs (i.e., Run 5 for both weighting matrices for Case E) are reported below as,
\begin{align*}
    \bm{K}_1 & = \begin{bmatrix} 6.5225 & 0 & 0 & 0 & 0 & 0 \\
                                 0 & 98.1494 & 0 & 0 & 0 & 0 \\
                                 0 & 0 & 8.9658 & 0 & 0 & 0 \\
                                 0 & 0 & 0 & 10.2037 & 0 & 0 \\ 
                                 0 & 0 & 0 & 0 & 97.3815 & 0 \\ 
                                 0 & 0 & 0 & 0 & 0 & 20.7142 \end{bmatrix}, \\
    \bm{K}_2 & = \begin{bmatrix} 39.4746 & 17.5941 & -2.0538 & -3.3242 & 0.1723 & -0.3125 \\
                                 17.5941 & 68.1822 & -3.7857 & -6.9744 & 0.1840 & -0.5589 \\
                                 -2.0538 & -3.7857 & 4.6930 & 0.5193 & 2.9905 & 0.2195 \\
                                 -3.3242 & -6.9744 & 0.5193 & 11.5512 & 0.4138 & -5.6421 \\
                                 0.1723 & 0.1840 & 2.9905 & 0.4138 & 77.1079 & 1.5671 \\
                                 -0.3125 & -0.5589 & 0.2195 & -5.6421 & 1.5671 & 81.3064 \end{bmatrix}.
\end{align*}

\section{Conclusion}
An efficient parameterization of a full positive-definite weighting matrix for Lyapunov control (LC) laws is presented. The parametrization is based on eigendecomposition of symmetric matrices, where the set of parameters consists of eigenvalues and hyperspherical angles used to construct orthonormal matrices.

Considering a full weighting matrix for LC laws offers a larger solution space by considering cross-coupling error terms. A variety of benchmark orbit transfer problems are solved with an LC law for the standard diagonal and the proposed full weighting matrix representations. Results showed that improved time-optimal solutions can be obtained with the full weighting matrix. The improvement is especially evident for classes of more complex and long-duration orbit transfer maneuvers.

One of our future works is to investigate combining our novel parametrization with Q-law and assessing if any improvements could be gained in both time-optimal and fuel-optimal solutions (since the coasting mechanism of Q-law allows for fuel-optimal solutions). Another future work is to investigate non-constant time- and/or state-dependent parameters within the full weighting matrix to enlarge the solution space even more.

\bibliographystyle{AAS_publication}
\bibliography{low_thrust_references}

\end{document}